\documentclass[12pt]{article}
\usepackage{amssymb}
\usepackage{amsthm}
\usepackage{amsmath}
\usepackage{graphicx}
 \newtheorem{thm}{Theorem}[section]
 
 \newtheorem{lem}[thm]{Lemma}
 \newtheorem{prop}[thm]{Proposition}
 \theoremstyle{definition}
 \newtheorem{defn}[thm]{Definition}
 \newtheorem{rem}[thm]{Remark}
 \numberwithin{equation}{section}

\begin{document}
\title{Global Low Regularity Solutions of Quasi-linear Wave Equations}
\author{Yi Zhou\thanks{School of Mathematical Sciences, Fudan University,
Shanghai 200433, P. R. China; Key Laboratory of Mathematics for
Nonlinear Sciences (Fudan University), Ministry of Education, P.
R. China {\it Email: yizhou@fudan.ac.cn.}} \and Zhen Lei
\thanks{School of Mathematical Sciences, Fudan University,
Shanghai 200433, P. R. China; School of Mathematics and
Statistics, Northeast Normal University, Changchun 130024, P. R.
China; Current Address: Applied and Comput. Math, Caltech,
Pasadena, CA 91125; School of Mathematics and Statistics,
Northeast Normal University, Changchun 130024, P. R. China. {\it
Email: zhenlei@acm.caltech.edu; leizhn@yahoo.com.}}}
\date{\today}
\maketitle


\begin{abstract}
In this paper we prove the global existence and uniqueness of the
low regularity solutions to the Cauchy problem of quasi-linear
wave equations with radial symmetric initial data in three space
dimensions. The results are based on the end-point Strichartz
estimate together with the characteristic method.
\end{abstract}



\section{Introduction}
This paper studies the global existence and uniqueness of the low
regularity solutions to the Cauchy problem associated to the
following quasi-linear wave equations in $\mathbb{R}^{1 + 3}$:
\begin{equation}\label{a1}
u_{tt} - a^2(u)\Delta u = 0,\quad t \geq 0,\ x \in \mathbb{R}^3,
\end{equation}
where $u = u(t, x)$ is a real-valued scalar function of $(t, x)$,
$u_{tt} = \partial_t^2u$, $\Delta$ is the Laplace operator acting
on space variables $x$, and $a$ is a smooth real-valued function
with $a(0)
> 0$. Equation \eqref{a1} is imposed on the following initial data
\begin{equation}\label{a2}
u(0, x) = f(x), u_t(0, x) = g(x),\quad  x \in \mathbb{R}^3.
\end{equation}

Before going any further, we shall first recall a few historical
results on the global classical solutions of system
\eqref{a1}-\eqref{a2}. In \cite{Lindblad1}, Lindblad proved the
global existence of classical solutions to system
\eqref{a1}-\eqref{a2} for sufficiently small, radially symmetric
initial data decaying rapidly as $x \rightarrow \infty$ or with
compact support. Alinhac \cite{Alinhac} improved the global
existence result of Lindblad \cite{Lindblad1} by removing the
assumption that the initial data is radially symmetric. Recently,
Lindblad \cite{Lindblad2} reproved the global existence result for
more general quasi-linear wave equations by a different and
simpler method.

On the other hand, many authors studied the local well-posedness
of low regularity solutions of quasi-linear wave equations
\cite{BahouriChemin1, BahouriChemin2, Klainerman,
KlainermanRodnianski}. Recently, Smith and Tataru
\cite{SmithTataru} established the local well-posedness of
quasi-linear wave equations in three space dimensions with initial
data in $H^s(\mathbb{R}^3)\times H^{s-1}(\mathbb{R}^3)$ for $s>2$.
This result is sharp in view of the counter example of Lindblad
\cite{Lindblad3}. We point out here that in this paper
$H^s(\mathbb{R}^3)$ denotes the usual Sobolev space endowed with
the norm
\begin{equation}\nonumber
\|f\|_{H^s(\mathbb{R}^3)} = \big\|(1 +
|\xi|^2)^{\frac{s}{2}}\hat{f}\big\|_{L^2(\mathbb{R}^3)},
\end{equation}
where $\hat{f}$ denotes the Fourier transform of $f$. For integer
$s \geq 0$, the above norm is equivalent to
\begin{equation}\label{a3}
\|f\|_{H^s(\mathbb{R}^3)} = \sum_{j =
0}^s\|\nabla^jf\|_{L^2(\mathbb{R}^3)}.
\end{equation}

Our main result in this paper is that the Cauchy problem
\eqref{a1}-\eqref{a2} is globally well-posed in the Sobolev space
$H^2(\mathbb{R}^3)\times H^1(\mathbb{R}^3)$ provided that the
initial data $(f, g)$ are \textit{radially symmetric} and have
compact supports, and the $H^2(\mathbb{R}^3)\times
H^1(\mathbb{R}^3)$ norm of the initial data $(f, g)$ is small
enough.

\begin{rem}\label{rem10}
Throughout this paper, we mean $\|\nabla f\|_{H^1(\mathbb{R}^3)} +
\|g\|_{H^1(\mathbb{R}^3)}$ by saying the $H^2(\mathbb{R}^3)\times
H^1(\mathbb{R}^3)$ norm of a pair $(f, g)$ and $\|\nabla u(t,
\cdot)\|_{H^1(\mathbb{R}^3)} + \|u_t(t,
\cdot)\|_{H^1(\mathbb{R}^3)}$ by saying the
$H^2(\mathbb{R}^3)\times H^1(\mathbb{R}^3)$ norm of a function
$u(t, x)$ at time $t$.
\end{rem}

To state our result more precisely, we introduce the concept of
strong solutions.

\begin{defn}\label{def1}
We say $u(t, x)$  is a strong solution to the Cauchy problem
\eqref{a1}-\eqref{a2} on some time interval $[0, T_0)$, if for a
sequence of initial data $(f_n, g_n) \in H^3(\mathbb{R}^3)\times
H^2(\mathbb{R}^3)$ which tends to $(f, g)$ strongly in
$H^2(\mathbb{R}^3)\times H^1(\mathbb{R}^3)$ as $n \rightarrow
\infty$, there exists a sequence of solutions $u_n(t, x)$ in
$C\big([0, T_0); H^3(\mathbb{R}^3)\big)\times C^1\big([0, T_0);
H^2(\mathbb{R}^3)\big)$ which solves the quasi-linear wave
equation \eqref{a1} with the initial data
\begin{equation}\nonumber
u_n(0, x) = f_n(x),\quad \partial_tu_n(0, x) = g_n(x),
\end{equation}
and tends to $u$ weakly$\star$ in $L^\infty\big([0, T],
H^2(\mathbb{R}^3)\big)\times W^{1, \infty}\big([0, T],
H^1(\mathbb{R}^3)\big)$ for any $T < T_0$ as $n \rightarrow
\infty$.

If $T_0 = + \infty$, then we say is $u(t, x)$  a global strong
solution to the Cauchy problem \eqref{a1}-\eqref{a2}. If $(f, g)$
and $(f_n, g_n)$ are all radially symmetric functions, then we say
$u(t, x)$ a radially symmetric strong solution to the Cauchy
problem \eqref{a1}-\eqref{a2}.
\end{defn}

\begin{rem}\label{rem11}
The local existence of $H^3(\mathbb{R}^3)\times H^2(\mathbb{R}^3)$
solutions for quasi-linear wave equations is guaranteed by the
classical local existence theorem provided that the initial data
is in $H^3(\mathbb{R}^3)\times H^2(\mathbb{R}^3)$. For example,
see Hughes-Kato-Marsden \cite{Hughes-Kato-Marsden}.
\end{rem}

Now we state our main result.
\begin{thm}\label{thm11}
Assume that

$(H_1)$:\quad $f$ and $g$ are radially symmetric functions with
compact supports:
\begin{equation}\label{a4}
{\rm supp}\big\{(f, g)\big\} \subseteq \big\{x \in \mathbb{R}^3,
|x| \leq 1\big\}.
\end{equation}

$(H_2)$:\quad $(f, g) \in H^2(\mathbb{R}^3) \times
H^1(\mathbb{R}^3)$ with
\begin{equation}\label{a5}
\|\nabla f\|_{H^1(\mathbb{R}^2)} + \|g\|_{H^1(\mathbb{R}^2)} \leq
\epsilon
\end{equation}
for a positive constant $\epsilon$.

Then there exists a unique global strong solution
\begin{equation}\label{a6}
u \in L^\infty_{{\rm loc}}\big([0, \infty); H^2(\mathbb{R}^3)\big)
\cap W^{1, \infty}_{{\rm loc}}\big([0, \infty);
H^1(\mathbb{R}^3)\big)
\end{equation}
to the quasi-linear wave equation \eqref{a1} with initial data
\eqref{a2} provided that $\epsilon$ is small enough. Moreover,
there exists a positive constant $A$ and a small enough positive
constant $\theta$ such that the global strong solution $u(t, x)$
satisfies
\begin{eqnarray}\label{a7}
&&\|\nabla u(t, \cdot)\|_{H^1(\mathbb{R}^2)}^2 + \|u_t(t,
  \cdot)\|_{H^1(\mathbb{R}^2)}^2\\\nonumber
&&\leq A\Big(\|\nabla f\|_{H^1(\mathbb{R}^2)}^2 +
  \|g\|_{H^1(\mathbb{R}^2)}^2\Big)(1 + t)^{\theta}
\end{eqnarray}
for all time $t \geq 0$, where $A$ depends only on the function
$a$ and $\theta$ depends only on $a$ and $\epsilon$.
\end{thm}

There are a lot of results on the global well-posedness of low
regularity solutions to semi-linear wave equations \cite{KeelTao,
KeelTao2, KlainermanMachedon2, KlainermanMachedon3, RodnianskiTao,
tao1}, however, to the best of our knowledge, theorem \ref{thm11}
is the only result for global existence and uniqueness of low
regularity solutions to quasi-linear wave equations.

The rest of this paper is organized as follows. In section 2, we
illustrate our main ideas of proving Theorem \ref{thm11}. More
precisely, we will smooth the initial data $(f, g)$ in \eqref{a2}
such that the local existence of solutions to the quasi-linear
wave equations \eqref{a1} with the regularized initial data is
guaranteed by Remark \ref{rem11}. Then we state the main
\textit{a\ priori} estimates \eqref{b18}-\eqref{b21} for the local
classical solutions involving only the $H^2(\mathbb{R}^3) \times
H^1(\mathbb{R}^3)$ norm of the initial data. The global existence
and the uniform bound \eqref{a7} for strong solutions to the
Cauchy problem \eqref{a1}-\eqref{a2} is a direct consequence of
\eqref{b21} (see Remark \ref{rem21}). We discuss how much the
characteristics of the quasi-linear wave equations \eqref{a1}
deviate from that of the linear wave equation in section 3, and
some other useful properties related to characteristics are also
proved. In section 4, we use the characteristic method to prove
the \textit{a\ priori} estimates presented in section 2 by
assuming that the weighted end-point Strichartz estimate
\eqref{b21} is true. To prove \eqref{b21}, we invoke Klainerman's
vector fields and the generalized energy method to explore decay
properties of energy away from the characteristics in section 5.
Then we estimate Hardy-Littlewood Maximal functions to get
\eqref{b21} in section 6. Section 7 studies the local existence of
general strong solutions to the Cauchy problem
\eqref{a1}-\eqref{a2} in Sobolev space $H^2(\mathbb{R}^3) \times
H^1(\mathbb{R}^3)$ and their stability in Sobolev space
$H^1(\mathbb{R}^3) \times L^2(\mathbb{R}^3)$, which implies the
uniqueness of the global strong solution and completes the proof
of theorem \ref{thm11}.

\section{A priori estimates}

This section illustrates our main observations and ideas in the
proof of Theorem \ref{thm11}, and also lays some groundwork for
the latter sections.

First of all, for any given radially symmetric function
\begin{equation}\nonumber
\rho(|x|) \in C_0^\infty(\mathbb{R}^3),\quad \rho \geq 0,\quad
\int_{\mathbb{R}^3}\rho(x)dx = 1,\quad \rho \equiv 0\ {\rm for}\
|x| \geq 2,
\end{equation}
we define the mollification $(\mathcal{J}_nf)(x)$ of functions $f
\in L^p(\mathbb{R}^3)$, $1 \leq p \leq \infty$, by
\begin{equation}\label{b1}
(\mathcal{J}_nf)(x) = n^3\int_{\mathbb{R}^3}\rho\big(n(x -
y)\big)f(y)dy,\quad n = 1, 2, 3, \cdots
\end{equation}
For any rotation matrix $Q$, it is easy to see that
\begin{eqnarray}\nonumber
&&(\mathcal{J}_nf)(Qx) = n^3\int_{\mathbb{R}^3}\rho\big(n(Qx -
  y)\big)f(y)dy\\\nonumber
&&= n^3\int_{\mathbb{R}^3}\rho\big(n(Qx -
  Qy)\big)f(Qy)dy\\\nonumber
&&= n^3\int_{\mathbb{R}^3}\rho\big(n(x -
  y)\big)f(Qy)dy.
\end{eqnarray}
Thus, for any radially symmetric function $f(x) = f(|x|)$, the
mollified function is also radially symmetric:
$(\mathcal{J}_nf)(x) = (\mathcal{J}_nf)(|x|)$.

Our first step of proving Theorem \ref{thm11} is to regularize the
initial data $(f, g) \in H^2(\mathbb{R}^3) \times
H^1(\mathbb{R}^3)$ so that the mollified functions $(f_n, g_n) =
\big(\mathcal{J}_nf, \mathcal{J}_ng\big) \in H^3(\mathbb{R}^3)
\times H^2(\mathbb{R}^3)$ by using the mollification operator
defined in \eqref{b1}. By the local existence theorem (see Remark
\ref{rem11}), there admits a unique solution $u_n(t, x)$ in
$H^3(\mathbb{R}^3) \times H^2(\mathbb{R}^3)$ to the quasi-linear
wave equation \eqref{a1} with the initial data $(f_n, g_n)$. Thus,
to show that there exists a unique global strong solution to the
Cauchy problem \eqref{a1}-\eqref{a2}, we only need to achieve the
\textit{a\ priori} estimate for $u_n$ in $H^2(\mathbb{R}^3) \times
H^1(\mathbb{R}^3)$ uniformly with respect to $n$, which only
involves the $H^2(\mathbb{R}^3) \times H^1(\mathbb{R}^3)$ norm
(but not the $H^3(\mathbb{R}^3) \times H^2(\mathbb{R}^3)$ norm) of
the initial data $(f_n, g_n)$. Note that $(f_n, g_n)$ is also
radially symmetric. In what follows, we will drop the subscript
$n$ and intend to achieve \textit{a\ priori} $H^2(\mathbb{R}^3)
\times H^1(\mathbb{R}^3)$ estimate for classical solutions $u(t,
x)$ to the quasi-linear wave equations \eqref{a1} only using the
$H^2(\mathbb{R}^3) \times H^1(\mathbb{R}^3)$ norm of the initial
data $(f, g)$.

We turn to focus on the standard energy estimates for the
quasi-linear wave equations \eqref{a1}. For any integer $s \geq
1$, define the standard energy by
\begin{equation}\label{b2}
E_s\big[u(t)\big] = \frac{1}{2}\sum_{b = 0}^{s -
1}\int_{\mathbb{R}^3} \Big\{|\partial^bu_t(t, \cdot)|^2 +
a^2(u)|\nabla\partial^bu(t, \cdot)|^2\Big\}dx,
\end{equation}
where $\partial = (\partial_t, \nabla)$. By the standard energy
estimates and Nirenberg inequality, it is rather easy to get (see
\cite{SmithTataru} and the references therein)
\begin{equation}\nonumber
\frac{d}{dt}E_s\big[(u(t)\big] \leq C_0\|\partial u(t,
\cdot)\|_{L^\infty(\mathbb{R}^3)}E_s\big[u(t)\big],
\end{equation}
provided that $u$ is bounded and $a^2(u)$ is strictly greater than
0, where $C_0$ is a positive constant depending only on $a$ and
$\|u\|_{L^\infty(\mathbb{R}^3)}$. Consequently, it follows that
\begin{equation}\nonumber
E_s\big[u(T)\big] \leq
E_s\big[u(0)\big]\exp\Big\{C_0\int_0^T\|\partial u(t,
\cdot)\|_{L^\infty(\mathbb{R}^3)}dt\Big\}.
\end{equation}

If the inequality
\begin{equation}\label{b3}
\int_0^T\|\partial u(t, \cdot)\|_{L^\infty(\mathbb{R}^3)}dt \leq
\frac{1}{K}\ln(1 + T)
\end{equation}
is valid for all time $T \geq 0$ and a uniform (big) constant $K >
0$, then there holds the following \emph{a priori} energy estimate
\begin{equation}\label{b4}
E_2\big[u(T)\big] \leq E_2\big[u(0)\big](1 + T)^{\frac{C_0}{K}},
\end{equation}
which only involves the $H^2(\mathbb{R}^3) \times
H^1(\mathbb{R}^3)$ norm of the initial data $(f, g)$. The energy
$E_2\big[u(T)\big]$ admits a dynamic growth, but will be still
bounded for all time $T \geq 0$. By the counterexample given by
Lindblad \cite{Lindblad3}, the inequality \eqref{b3} may not be
true for the quasi-linear wave equation \eqref{a1} with general
initial data $(f, g)$ in $H^2(\mathbb{R}^3) \times
H^1(\mathbb{R}^3)$. We will prove in this paper that the
inequality \eqref{b3} is true (see Remark \ref{rem21}) for all
time $t \geq 0$ if the initial data $(f, g)$ satisfies the
assumptions $(H_1)-(H_2)$ in Theorem \ref{thm11}. The global
existence and the uniform bound \eqref{a7} for strong solutions to
the Cauchy problem of quasi-linear wave equation
\eqref{a1}-\eqref{a2} are then direct consequences of \eqref{b4}
(by letting $\theta = \frac{C_0}{K}$), as what we have mentioned
at the beginning of this section.

To present our main ideas of proving the \textit{a\ priori}
estimate \eqref{b3}, we begin with defining the characteristics of
the quasi-linear wave equations \eqref{a1}.

Let $r_-(\tau; \beta)$ denote a family of minus characteristics
passing through the point $(t, r)$ with the intersection $(0,
\beta)$, $\beta \geq 0$, with $\sigma$ axis (see Figure 1):
\begin{equation}\label{b5}
\begin{cases}
\frac{dr_-(\tau; \beta)}{d\tau} = - a(u)\big(\tau, r_-(\tau;
  \beta)\big),\\[-3mm]\\
r_-(0; \beta) = \beta.
\end{cases}
\end{equation}
For convenience, sometimes we also denote the above minus
characteristics as $t_-(\sigma; \beta)$  (see Figure 1):
\begin{equation}\label{b6}
\begin{cases}
\frac{dt_-(\sigma; \beta)}{d\sigma} = - \frac{1}
  {a(u)\big(t_-(\sigma; \beta), \sigma\big)},\\[-4mm]\\
t_-(\beta; \beta) = 0.
\end{cases}
\end{equation}

\begin{figure}
\medskip
\includegraphics [width=12cm,clip]{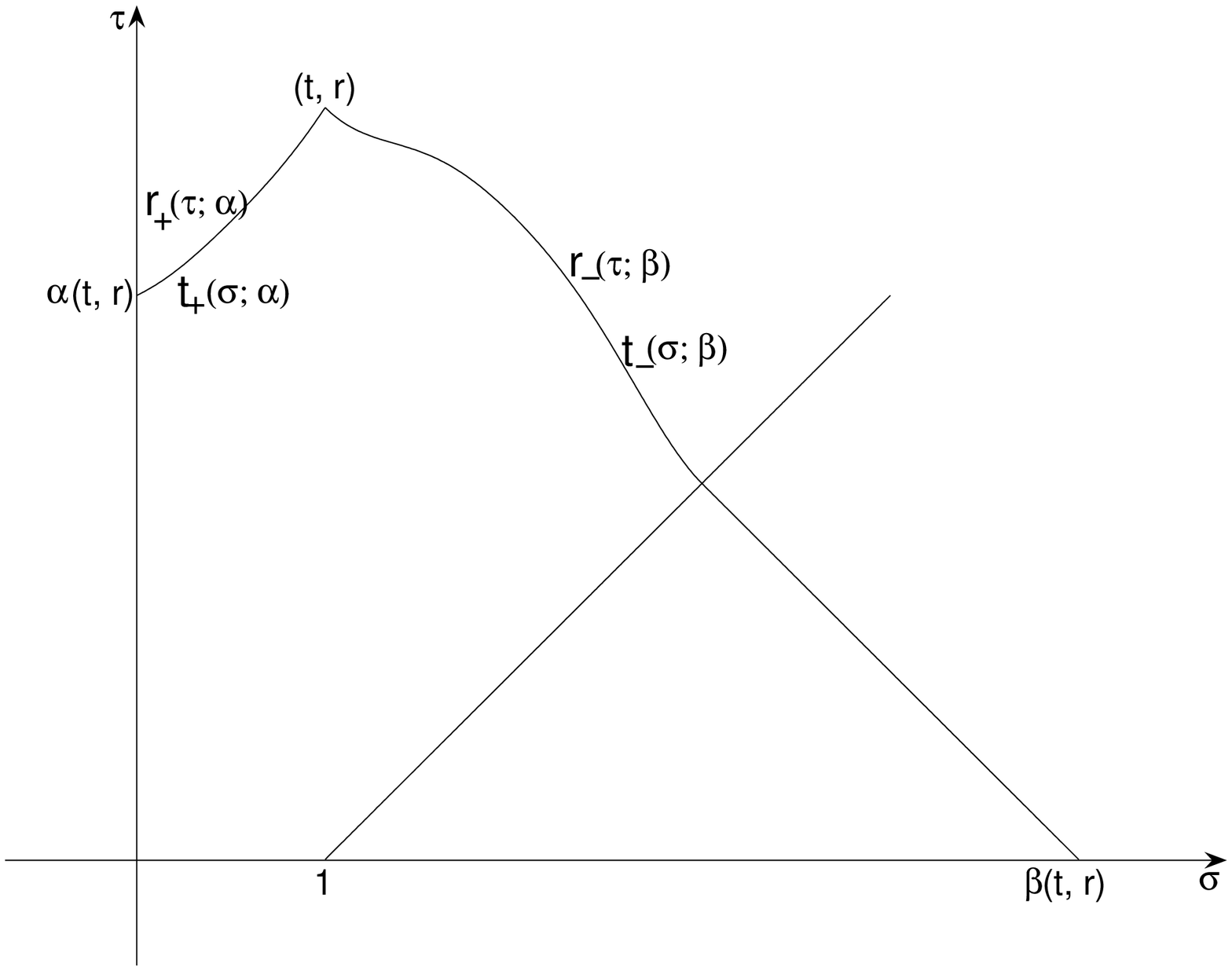}
\caption{}
\medskip
\end{figure}

Similarly, let $r_+(\tau; \alpha)$ denote a family of plus
characteristics passing through the point $(t, r)$ with the
intersection $(\alpha, 0)$, $\alpha \geq 0$, with $\tau$ axis (see
Figure 1):
\begin{equation}\label{b7}
\begin{cases}
\frac{dr_+(\tau; \alpha)}{d\tau} = a(u)\big(\tau,
r_+(\tau; \alpha)\big),\\[-3mm]\\
r_+(\alpha; \alpha) = 0.
\end{cases}
\end{equation}
which is also written as $t_+(\sigma; \alpha)$ (see Figure 1):
\begin{equation}\label{b8}
\begin{cases}
\frac{dt_+(\sigma; \alpha)}{d\sigma} = \frac{1}
{a(u)\big(t_+(\sigma; \alpha), \sigma\big)},\\[-3mm]\\
t_+(0; \alpha) = \alpha,
\end{cases}
\end{equation}
When starting from $(0, \gamma)$, $\gamma > 0$, a plus
characteristic $r_+(\tau; \gamma)$ is defined by
\begin{equation}\label{b9}
\begin{cases}
\frac{dr_+(\tau; \gamma)}{d\tau} = a(u)\big(\tau,
r_+(\tau; \gamma)\big),\\[-3mm]\\
r_+(0; \gamma) = \gamma,
\end{cases}
\end{equation}
and sometimes we also write it as $t_+(\sigma; \gamma)$:
\begin{equation}\label{b10}
\begin{cases}
\frac{dt_+(\sigma; \gamma)}{d\sigma} = \frac{1}
{a(u)\big(t_+(\sigma; \gamma), \sigma\big)},\\[-3mm]\\
t_+(\gamma; \gamma) = 0.
\end{cases}
\end{equation}

Note that in the above definitions \eqref{b5}-\eqref{b10}, the
parameters $\alpha$, $\beta$ and $\gamma$ depend on the point $(t,
r)$, and the point $(t, r)$ also depends on the parameters
$(\alpha, \beta)$ or $(\gamma, \beta)$.  To show their dependence,
we also write $\alpha = \alpha(t, r)$, $\beta = \beta(t, r)$,
$\gamma = \gamma(t, r)$ and $t = t(\alpha, \beta)$, $r = r(\alpha,
\beta)$ or $t = t(\gamma, \beta)$, $r = r(\gamma, \beta)$. It is
obviously that
\begin{equation}\label{b11}
\begin{cases}
t_+\big(r; \alpha(t, r)\big) = t,\quad r_+\big(t; \alpha(t,
r)\big) =
r,\\[-4mm]\\
t_+\big(r; \gamma(t, r)\big) = t,\quad r_+\big(t; \gamma(t, r)\big) = r,\\[-4mm]\\
t_-\big(r; \beta(t, r)\big) = t,\quad r_-\big(t; \beta(t, r)\big) = r,\\[-4mm]\\
\alpha\big(\tau, r_+(\tau; \alpha)\big) = \alpha\big(t_+(\sigma;
\alpha), \sigma\big) = \alpha,\\[-4mm]\\
\gamma\big(\tau, r_+(\tau; \gamma)\big) = \gamma\big(t_+(\sigma;
\gamma), \sigma\big) = \gamma,\\[-4mm]\\
\beta\big(\tau, r_-(\tau; \beta)\big) = \beta\big(t_-(\sigma;
\beta), \sigma\big) = \beta.
\end{cases}
\end{equation}
These notations seem rather tedious and complicated. However, the
reader will see the advantages of the above convention.

Next, we define the cone
\begin{equation}\nonumber
\mathcal{D} = \{(t, r)|r \leq \frac{t}{4} + 1\},
\end{equation}
and assume that (without loss of generality)
\begin{equation}\label{b16}
a(0) = 1.
\end{equation}
As in \cite{Lindblad1}, we introduce the following weighted
differential operators:
\begin{equation}\label{b12}
L_\pm = \frac{1}{\sqrt{a(u)}}\partial_t \pm \sqrt{a(u)}\partial_r.
\end{equation}

We will prove the following crucial Theorem:
\begin{thm}\label{thm21}
Suppose that the assumptions $(H_1)-(H_2)$ in Theorem \ref{thm11}
are satisfied and $u(t, r)$ is the classical solution to the
Cauchy problem of the quasi-linear wave equation
\eqref{a1}-\eqref{a2}. Suppose furthermore that $K$ ($\gg$ 1) and
$\epsilon$ ($\ll$ 1) satisfy
\begin{equation}\label{b17}
K^4\epsilon \leq 1.
\end{equation}
Then there exists positive constants $C_j> 0$, $j = 1, 2, \cdots,
6$, such that for all time $t \geq 0$, there holds

A priori estimates for $u$:
\begin{equation}\label{b18}
\begin{cases}
|u(t, r)| \leq \frac{C_1\epsilon}{(1 + t)^{\frac{3}{5}}},\\[-3mm]\\
|u(t, r)| \leq C_2K\epsilon(1 + t)^{- 1 + \frac{2}{K}}.
\end{cases}
\end{equation}

A priori estimates for $L_\pm v$ with $v = ru$:
\begin{equation}\label{b19}
\begin{cases}
|L_+v(t, r)| \leq \frac{C_3\epsilon(1 + t)^{\frac{1}{K}}}
  {1 + \beta(t, r)},\\[-3mm]\\
|L_-v(t, r)| \leq \frac{C_4\epsilon}{\big(1 + \gamma(t, r)
  \big)^{1 - \frac{1}{K}}} \leq C_4\epsilon \quad {\rm
  for}\quad t_+(r; 1) \leq t < t_+(r; 0),\\[-3mm]\\
|L_-v(t, r)| \leq \frac{C_4\epsilon}{\big(1 + \alpha(t, r)\big)^{1
  - \frac{1}{K}}} \quad {\rm for}\quad t \geq t_+(r; 0).
\end{cases}
\end{equation}

A priori estimates for $\partial u$ holds away from the cone
$\mathcal{D}$:
\begin{equation}\label{b20}
\begin{cases}
|\partial u(t, r)| \leq \frac{C_5K\epsilon}{(1 + t)\big(1 +
  \gamma(t, r)\big)^{1 - \frac{2}{K}}},\quad{\rm for}\quad
  r > \max\{\frac{t}{4} + 1, r_+(t; 0)\},\\[-3mm]\\
|\partial u(t, r)| \leq \frac{C_5K\epsilon}{(1 + t)\big(1 +
  \alpha(t, r)\big)^{1 - \frac{2}{K}}},\quad{\rm for}\quad
  \frac{t}{4} + 1 < r \leq r_+(t; 0).
\end{cases}
\end{equation}

Weighted a priori estimate for $\partial u$ holds inside the cone
$\mathcal{D}$:
\begin{equation}\label{b21}
\int_0^t\big[(1 + \tau)^{1 - \frac{4}{K}}\sup_{\sigma \leq
\frac{\tau}{4} + 1}|\partial u(\tau, \sigma)|\big]^2d\tau \leq
C_6^2K^4\epsilon^2.
\end{equation}
\end{thm}

\begin{rem}\label{rem21}
With the aid of \eqref{b17}, estimates \eqref{b20} and \eqref{b21}
indeed imply \eqref{b3} by the following calculations:
\begin{eqnarray}\nonumber
&&\int_0^t|\partial u(\tau, \cdot)|_{L^\infty} d\tau\\\nonumber
&&\leq \int_0^t\big[\sup_{\sigma >
  \frac{\tau}{4} + 1}|\partial u(\tau, \sigma)|\big]
  d\tau + \int_0^t\big[\sup_{\sigma \leq
  \frac{\tau}{4} + 1}|\partial u(\tau, \sigma)|\big]
  d\tau\\\nonumber
&&\leq C_5K\epsilon\ln(1 + t) + C\big[1 - (1 + t)^{-
   1 + \frac{8}{K}}\big]^{\frac{1}{2}}\\\nonumber
&&\quad \times \Big(\int_0^t\big[(1 + \tau)^{1
  - \frac{4}{K}}\sup_{\sigma \leq \frac{\tau}{4} + 1}|\partial
  u(\tau, \sigma)|\big]^2d\tau\Big)^{\frac{1}{2}}\\\nonumber
&&\leq CK\epsilon\ln(1 + t) + CK^2\epsilon\big[1 - (1 + t)^{-
   1 + \frac{8}{K}}\big]^{\frac{1}{2}}\\\nonumber
&&\leq CK^2\epsilon\ln(1 + t).
\end{eqnarray}
\end{rem}

We mention that in this paper, $C$ will be used to denote a
generic positive constant independent of any function appeared
throughout this paper except for $a$, and its meaning may vary
from line to line.

The proof of \eqref{b18} will be completed by using \eqref{b19} at
the beginning of section 4. \eqref{b19} will be proved by using
\eqref{b21} and the characteristic method also in section 4. The
properties of characteristics will be discussed in section 3. Then
we use \eqref{b18} and \eqref{b19} to prove \eqref{b20} at the end
of section 4. The strategy of proving the series of estimates
\eqref{b18}-\eqref{b21} is that by assuming their validity first,
we show these estimates are still valid when the constants $C_j$
are replaced by $\frac{1}{2}C_j$ for $j = 1, 2, \cdots, 6$.

The proof of \eqref{b21} relies on the estimation of
Hardy-Littlewood Maximal functions and involves Klainerman's
vector fields and generalized energy method. To best illustrate
our ideas, let us review the end-point Strichartz estimate by
Klainerman and Machedon \cite{KlainermanMachedon}. Consider the
following linear problem
\begin{equation}\label{b22}
\begin{cases}
\phi_{tt} - \Delta \phi = 0,\quad t \geq 0, x \in
\mathbb{R}^3,\\[-4mm]\\
\phi(0, x) = \phi_0(r), \phi_t(0, x) = 0,\quad x \in \mathbb{R}^3.
\end{cases}
\end{equation}

Extend the initial data $\phi_0(r)$ and the solution $\phi(t, r)$
of system \eqref{b22} into even functions with respect to $r$, and
then let
\begin{equation}\label{b23}
\widetilde{\phi}_0 = r\phi_0(r), \quad r \in \mathbb{R}^1,.
\end{equation}
d'Alembert formula gives
\begin{eqnarray}\nonumber
&&\phi(t, r)\\\nonumber &&= \frac{\widetilde{\phi}_0(r + t) +
  \widetilde{\phi}_0(r - t)}{2r}\\\nonumber
&&= \frac{\widetilde{\phi}_0(t + r) -
  \widetilde{\phi}_0(t - r)}{2r}\\\nonumber
&&= \frac{1}{2r}\int_{t - r}^{t + r}\widetilde{\phi}_0
  ^\prime(\lambda)d\lambda.
\end{eqnarray}
Differentiating the above equality with respect to $t$ yields
\begin{equation}\label{b24}
\phi_t(t, r) = \frac{1}{2r}\int_{t - r}^{t +
r}\widetilde{\phi}_0^{\prime\prime}(\lambda)d\lambda.
\end{equation}

We recall some properties for Hardy-Littlewood Maximal functions.
The Hardy-Littlewood maximal operator in $\mathbb{R}^n$ is defined
on $L^1_{{\rm loc}}(\mathbb{R}^n)$ by
\begin{equation}\nonumber
\mathbb{M}f(x) = \sup_Q\frac{1}{|Q|}\int_Q|f(y)|dy,\quad f \in
L^1_{{\rm loc}}(\mathbb{R}^n),
\end{equation}
where the supremum is taken over all cubes containing $x$. The
following property is well-known (see \cite{stein}):
\begin{prop}\label{prop1}
The Hardy-Littlewood maximal operator is of strong $(p, p)$ type
for $1 < p \leq \infty$:
\begin{equation}\nonumber
\mathbb{M}f \in L^p(\mathbb{R}^n),\
\|\mathbb{M}f\|_{L^p(\mathbb{R}^n)} \leq
C\|f\|_{L^p(\mathbb{R}^n)}\quad {\rm for}\quad f \in
L^p(\mathbb{R}^n), 1 < p \leq \infty,
\end{equation}
where $C > 0$ is a uniform constant.
\end{prop}

To simplify our notations, we will denote the $L^p$ and Sobolev
norms of $\phi(t, r)$ with respect to $r$ by $\|\phi(t,
\cdot)\|_{L^p}$ and $\|\phi(t, \cdot)\|_{H^s}$, compared with the
notations for the Sobolev norms of a function with respect to $x
\in \mathbb{R}^3$ which are denoted by $\|\phi(t,
\cdot)\|_{L^p(\mathbb{R}^3)}$ and $\|f(t,
\cdot)\|_{H^s(\mathbb{R}^3)}$ (see \eqref{a3}).

Using \eqref{b24}, Hardy inequality and Proposition \ref{prop1},
we have (see Klainerman and Machedon \cite{KlainermanMachedon})
\begin{eqnarray}\label{b25}
&&\Big(\int_0^T\|\phi_t(t, \cdot)\|_{L^\infty}^2
  dt\Big)^{\frac{1}{2}}\\\nonumber
&&\leq \Big(\int_0^T\big|\mathbb{M}\big[\widetilde{\phi}_0^{
  \prime\prime}(\lambda)\big](t)\big|^2dt\Big)^{\frac{1}{2}}\\\nonumber
&&\leq C\Big(\int_0^\infty\big[|\phi_0^{
  \prime}(\lambda)|^2 + |\lambda\phi_0^{
  \prime\prime}(\lambda)|^2\big]d\lambda\Big)^{\frac{1}{2}}\\\nonumber
&&\leq C\|\phi_0^{\prime\prime}\|_{L^2(\mathbb{R}^3)} +
  C\big\|\frac{\phi_0^\prime(x)}{|x|}\big\|_{L^2(\mathbb{R}^3)}\\\nonumber
&&\leq C\|\phi_0^{\prime\prime}\|_{L^2(\mathbb{R}^3)}.
\end{eqnarray}
One of our key observations is that there holds weighted end-point
Strichartz estimate inside the cone $\mathcal{D}$:
\begin{eqnarray}\label{b26}
&&\Big(\int_0^T\big[(1 + t)^\delta\|\phi_t(t, r)\|_{L^\infty(r
  \leq \frac{t}{4} + 1)}\big]^2dt\Big)^{\frac{1}{2}}\\\nonumber
&&\leq \Big(\int_0^T\big|\mathbb{M}\big[(1 +
  \lambda)^\delta\widetilde{\phi}_0^{\prime\prime}(\lambda)\big]
  (t)\big|^2dt\Big)^{\frac{1}{2}}\\\nonumber
&&\leq C\Big(\int_0^\infty\big[|\phi_0^{
  \prime}(\lambda)|^2 + |\lambda\phi_0^{
  \prime\prime}(\lambda)|^2\big]d\lambda\Big)^{\frac{1}{2}}\\\nonumber
&&\leq C\|\phi_0^{\prime\prime}\|_{L^2(\mathbb{R}^3)}
\end{eqnarray}
provided that $\phi_0$ has a compact support, where $\delta$ is a
positive constant. Using equality \eqref{b24} and estimate
\eqref{b26}, we have
\begin{eqnarray}\label{b29}
&&\int_0^T\|\phi_t(t, r)\|_{L^\infty}dt\\\nonumber &&\leq
  \int_0^T\|\phi_t(t, r)\|_{L^\infty(r \leq \frac{t}{4} + 1)}dt
  + \int_0^T\|\phi_t(t, r)\|_{L^\infty(r
  > \frac{t}{4} + 1)}dt\\\nonumber
&&\leq \frac{C\big[1 - (1 + t)^{1 - 2\delta}\big]^{
  \frac{1}{2}}}{2\delta - 1}\Big(\int_0^T\big[(1 +
  t)^\delta\|\phi_t(t, r)\|_{L^\infty(r \leq \frac{t}{4}
  + 1)}\big]^2dt\Big)^{\frac{1}{2}}\\\nonumber
&&\quad +\ C\|\phi_0^{\prime\prime}\|_{L^1}\ln(1 + t)\\\nonumber
&&\leq C(\delta)\|\phi_0^{\prime\prime}\|_{L^2(\mathbb{R}^3)}\ln(1
+ t)
\end{eqnarray}
for $\frac{1}{2} < \delta < 1$.

To obtain the \textit{a priori} estimate \eqref{b21} for
quasi-linear wave equation \eqref{a1}, we need inevitably estimate
the second order derivatives of $u$, as has been partly revealed
in \eqref{b24}. This requires to invoke the Klainerman's vector
fields and generalized energy method to explore the decay
properties of energy inside the cone $\mathcal{D}$, which will be
discussed in section 5.

To use the generalized energy to estimate the maximal function, we
need introduce the following special weighted differential
operators
\begin{eqnarray}\label{b27}
M_{\pm} = \frac{1}{a(u)}\partial_t \pm \partial_r,
\end{eqnarray}
and rewrite the equation in \eqref{a1} as (see the derivation in
\eqref{f1}-\eqref{f2})
\begin{eqnarray}\label{b28}
&&M_-M_+[(t - r) v]\\\nonumber&&= - \frac{a^\prime(u)(t -
r)rM_-uu_t}{a^2(u)} + \Big(\frac{1 +
  a(u)}{a(u)}\Big)M_+v + M_-\Big(\frac{\big(1
  - a(u)\big)v}{a(u)}\Big).
\end{eqnarray}
We point out that the term $r$ in the first term of the right hand
side of the above special form of equation \eqref{a1} is essential
so that one can use energy to estimate the Hardy-Littlewood
Maximal functions. The weight $t - r$ is also a key ingredient
which provides us extra decay for $\partial u$ inside the cone
$\mathcal{D}$ by using the fact that the accumulation of the
quantity $(1 + \alpha)^{- 2 + \frac{10}{K}}$ along a fixed minus
characteristic is not so big (see Lemma \ref{lem35}). These
combining the fact that energy is decaying faster inside the cone
$\mathcal{D}$ explored by the generalized energy estimate finally
yield the desired weighted estimate \eqref{b21}. We will discuss
the generalized energy method in section 5, and prove the weighted
estimate \eqref{b21} in section 6.

\section{Properties of Characteristics}

In section 2, we defined the families of plus and minus
characteristics in \eqref{b5}-\eqref{b10} which satisfy the basic
identities \eqref{b11}. We also stated that the classical solution
to the Cauchy problem of the quasi-linear wave equation
\eqref{a1}-\eqref{a2} satisfies the \textit{a priori} estimates in
Theorem \ref{thm21}. In this section, we use the \textit{a priori}
estimates of the classical solution claimed in Theorem \ref{thm21}
to study the useful properties of characteristics defined in
\eqref{b5}-\eqref{b10}.

The following Lemma will be used to estimate the derivatives of
characteristics with respect to the parameters $\alpha$, $\beta$
and $\gamma$.

\begin{lem}\label{lem31}
Let $u$ be a classical solution of the Cauchy problem
\eqref{a1}-\eqref{a2}, and $r_-(\tau; \beta)$, $t_-(\sigma;
\beta)$, $r_+(\tau; \alpha)$, $t_+(\sigma; \alpha)$, $r_+(\tau;
\gamma)$ and $t_+(\sigma; \gamma)$ be defined by
\eqref{b5}-\eqref{b10}. Then we have
\begin{equation}\label{c1}
\begin{cases}
r_{+, \alpha}(\tau; \alpha) = - a\big(u(\alpha,
  0)\big)\exp\Big\{\int_{\alpha}^\tau a^\prime(u)
  u_r\big(s, r_+(s; \alpha)\big)ds\Big\},\\[-4mm]\\
t_{+, \alpha}(\sigma; \alpha) = \exp\Big\{- \int_0^\sigma
  \frac{a^\prime(u)}{a^2(u)}u_t\big(t_+(\rho; \alpha),
  \rho\big)d\rho\Big\},\\[-4mm]\\
r_{+, \gamma}(\tau; \gamma) = \exp\Big\{\int_0^\tau a^\prime(u)
  u_r\big(s, r_+(s; \gamma)\big)ds\Big\},\\[-4mm]\\
t_{+, \gamma}(\sigma; \gamma) = -\frac{1}{a\big(f(
  \gamma)\big)}\exp\Big\{- \int_\gamma^\sigma
  \frac{a^\prime(u)}{a^2(u)}u_t\big(t_+(\rho; \gamma),
  \rho\big)d\rho\Big\},\\[-4mm]\\
r_{-, \beta}(\tau; \beta) = \exp\Big\{-\int_0^\tau
  a^\prime(u)u_r\big(s, r_-(s; \beta)\big)ds\Big\},\\[-4mm]\\
t_{-, \beta}(\sigma; \beta) = \frac{1}{a\big(f(
  \beta)\big)}\exp\Big\{\int_\beta^\sigma \frac{a^\prime(u)}{a^2(u)}
  u_t\big(t_-(\rho; \beta), \rho\big)d\rho\Big\}.
\end{cases}
\end{equation}
\end{lem}
\begin{proof}
These equalities follow by differentiating the characteristics in
\eqref{b5}-\eqref{b10} with respect to the parameters $\alpha$,
$\beta$ and $\gamma$, and then solving the resulting ordinary
differential equations along the corresponding plus or minus
characteristics.
\end{proof}

\begin{lem}\label{lem32}
Let $u$ be a classical solution of the Cauchy problem
\eqref{a1}-\eqref{a2}, and $r_-(\tau; \beta)$, $r_+(\tau; \alpha)$
and $r_+(\tau; \gamma)$ be defined in \eqref{b5}, \eqref{b7} and
\eqref{b9}. Then we have
\begin{equation}\label{c2}
\begin{cases}
\alpha_t(t, r) = - \frac{a\big(u(t, r)\big)}{r_{+, \alpha}\big(t;
\alpha(t, r)\big)},\quad \alpha_r(t, r) = \frac{1}{r_{+,
\alpha}\big(t; \alpha(t, r)\big)},\\[-3mm]\\
\gamma_t(t, r) = - \frac{a\big(u(t, r)\big)}{r_{+, \gamma}\big(t;
\gamma(t, r)\big)},\quad \gamma_r(t, r) = \frac{1}{r_{+,
\gamma}\big(t; \gamma(t, r)\big)},\\[-3mm]\\
\beta_t(t, r) = \frac{a\big(u(t, r)\big)}{r_{-, \beta}\big(t;
\beta(t, r)\big)},\quad \beta_r(t, r) = \frac{1}{r_{-,
\beta}\big(t; \beta(t, r)\big)},
\end{cases}
\end{equation}
and
\begin{equation}\label{c3}
\begin{cases}
t_\alpha(\alpha, \beta) = - \frac{r_{+, \alpha}\big(t(\alpha,
\beta); \alpha\big)}{2a(u)\big(t(\alpha, \beta), r(\alpha,
\beta)\big)},\quad t_\beta(\alpha, \beta) = \frac{r_{-,
\beta}\big(t(\alpha, \beta); \beta\big)} {2a(u)\big(t(\alpha,
\beta), r(\alpha, \beta)\big)},\\[-3mm]\\
r_\alpha(\alpha, \beta) = \frac{r_{+, \alpha}\big(t(\alpha,
\beta); \alpha\big)}{2},\quad r_\beta(\alpha, \beta) = \frac{r_{-,
\beta}\big(t(\alpha, \beta); \beta\big)}{2},\\[-3mm]\\
t_\gamma(\gamma, \beta) = - \frac{r_{+, \gamma}\big(t(\gamma,
\beta); \gamma\big)}{2a(u)\big(t(\gamma, \beta), r(\gamma,
\beta)\big)},\quad r_\gamma(\gamma, \beta) = \frac{r_{+,
\gamma}\big(t(\gamma, \beta); \gamma\big)}{2}.
\end{cases}
\end{equation}
\end{lem}
\begin{proof}
The identities in \eqref{c2} are the direct consequences of
differentiating the following identities (see \eqref{b11})
\begin{equation}\nonumber
r_+\big(t; \alpha(t, r)\big) = r,\quad r_+\big(t; \gamma(t,
r)\big) = r,\quad r_-\big(t; \beta(t, r)\big) = r
\end{equation}
with respect to $t$ and $r$. The identities in \eqref{c3} follow
by the following calculations:
\begin{equation}\nonumber
\frac{\partial (t, r)}{\partial(\alpha, \beta)} =
\begin{pmatrix}\alpha_t & \alpha_r\\ \beta_t &
\beta_r\end{pmatrix}^{- 1} =
\frac{1}{\alpha_t\beta_r - \beta_t\alpha_r}\begin{pmatrix}\beta_r & - \alpha_r\\
- \beta_t & \alpha_t\end{pmatrix},
\end{equation}
\end{proof}

The following Lemma tells us how much the characteristics of the
quasi-linear wave equation \eqref{a1} deviate from that of the
corresponding linear wave equation (see Figure 1).

\begin{lem}\label{lem33}
Suppose that \eqref{b17} and the assumptions $(H_1)-(H_2)$ in
Theorem \ref{thm11} are satisfied. Then we have
\begin{equation}\label{c4}
\begin{cases}
|r_+^\prime(\tau; \alpha) - 1| \leq C\epsilon,\quad
|r_-^\prime(\tau; \beta) + 1| \leq
C\epsilon,\\[-3mm]\\
|t_+^\prime(\sigma; \alpha) - 1| \leq C\epsilon,\quad
|t_-^\prime(\sigma; \beta) + 1| \leq C\epsilon,\\[-3mm]\\
|r_+^\prime(\tau; \gamma) - 1| \leq C\epsilon,\quad
|t_+^\prime(\sigma; \gamma) - 1| \leq C\epsilon,
\end{cases}
\end{equation}
and
\begin{equation}\label{c5}
\begin{cases}
\big|(t + r) - \beta(t, r)\big| \leq C\epsilon\min\{t,
K^2(1 + t)^{\frac{2}{K}}, \beta(t, r) - r\},\\[-3mm]\\
\big|(t - r) - \alpha(t, r)\big| \leq C\epsilon\min\{r,
K^2(1 + t)^{\frac{2}{K}}, t - \alpha(t, r)\},\\[-3mm]\\
\big|(r - t) - \gamma(t, r)\big| \leq C\epsilon\min\{t, K^2(1 +
t)^{\frac{2}{K}}, r - \gamma(t, r)\}.
\end{cases}
\end{equation}
\end{lem}

\begin{proof}
The first line of \eqref{c4} follows from the first inequality of
\eqref{b18} and the following calculation
\begin{eqnarray}\nonumber
&&|r_+^\prime(\tau; \alpha) - 1| + |r_-^\prime(\tau; \beta) +
  1|\\\nonumber
&&= \big|a\big[u\big(\tau, r_+(\tau; \alpha)\big)\big] - a(0)\big|
+
  \big|a\big[u\big(\tau, r_-(\tau;
  \beta)\big)\big] - a(0)\big|\\\nonumber
&&\leq C\big[|u\big(\tau, r_+(\tau; \alpha)\big)| + |u\big(\tau,
r_-(\tau;
  \beta)\big)|\big]\\\nonumber
&&\leq C\epsilon.
\end{eqnarray}
Similarly, one can prove other inequalities of \eqref{c4}.

To show the first estimate in \eqref{c5}, we compute
\begin{eqnarray}\nonumber
&&\big|t - \big(\beta(t, r) - r\big)\big|\\\nonumber &&=
  \big|\int_r^{\beta(t, r)}\Big\{\frac{1}{a(u)\big[t_-\big(\sigma;
  \beta(t, r)\big), \sigma\big]} - 1\Big\}d\sigma\big|\\\nonumber
&&= \big|\int_0^t\Big\{1 - a(u)\big[\tau, r_-\big(\tau;
  \beta(t, r)\big)\big]\Big\}d\tau\big|\\\nonumber
&&\leq C\epsilon\min\{t, \beta(t, r) - r\},
\end{eqnarray}
where we used the first inequality of \eqref{b18} and $a(0) = 1$.
Similarly, by using the second inequality instead of the first one
of \eqref{b18}, we also have
\begin{eqnarray}\nonumber
&&\big|t - \big(\beta(t, r) - r\big)\big|\\\nonumber &&\leq
CK\epsilon\big|\int_0^t(1 + \tau)^{- 1 +
\frac{2}{K}}d\tau\big|\\\nonumber &&\leq CK^2\epsilon(1 +
t)^{\frac{2}{K}}.
\end{eqnarray}
The proofs of the second and third inequalities of \eqref{c5} are
similar as above.
\end{proof}

By \eqref{c4}, we can easily get
\begin{lem}\label{lem34}
Suppose that \eqref{b17} and the assumptions $(H_1)-(H_2)$ in
Theorem \ref{thm11} are satisfied. Then for any $(t, r)$ with $t
\geq t_+(r; 1)$,
\begin{equation}\label{c6}
\begin{cases}
\frac{t + 1}{2} \leq 1 + \beta(t, r) \leq 3(t + 1),\\[-3mm]\\
\frac{t + 1}{2} \leq \alpha(t, r) + 1 \leq t + 1\quad {\rm
for}\quad
  r \leq \frac{t}{4} + 1,\\[-3mm]\\
\tau \leq t \leq 3\tau \quad {\rm for}\quad
  \tau \geq t_1,
\end{cases}
\end{equation}
where $t_1$ is determined by (see Figure 2)
\begin{equation}\label{c7}
r_+(t_1; 1) = r_-\big(t_1; \beta(t, r)\big).
\end{equation}
\end{lem}

The following Lemma describes the accumulation of the quantities
$(1 + \alpha)^{- 2 + \frac{10}{K}}$ (or $(1 + \gamma)^{- 2 +
\frac{10}{K}}$) and $(1 + \beta)^{- 2 + \frac{10}{K}}$ along the
minus and plus characteristics.

\begin{lem}\label{lem35}
Suppose that \eqref{b17} and the assumptions $(H_1)-(H_2)$ in
Theorem \ref{thm11} are satisfied. Then we have
\begin{equation}\label{c8}
\begin{cases}
\int_\alpha^t\frac{1}{\big[1 + \beta\big(\tau, r_+(\tau;
\alpha)\big)\big]^{2 - \frac{10}{K}}}d\tau \leq C,
\\[-3mm]\\
\int_0^t\frac{1}{\big[1 + \beta\big(\tau, r_+(\tau;
\gamma)\big)\big]^{2 - \frac{10}{K}}}d\tau \leq C,
\\[-3mm]\\
\int_{t_1}^t\frac{1}{\big[1 + \gamma\big(\tau, r_-(\tau;
\beta)\big)\big]^{2 - \frac{10}{K}}}d\tau \leq C(1 +
t)^{CK^2\epsilon}\quad {\rm for}\ t \leq t_2,\\[-3mm]\\
\int_{t_2}^t\frac{1}{\big[1 + \alpha\big(\tau, r_-(\tau;
\beta)\big)\big]^{2 - \frac{10}{K}}}d\tau \leq C(1 +
t)^{CK^2\epsilon} \quad {\rm for}\ t \geq t_2,
\end{cases}
\end{equation}
where $t_1$ is determined in \eqref{c7} and $t_2$ is defined by
(see Figure 2)
\begin{equation}\label{c9}
r_+(t_2; 0) = r_-(t_2, \beta).
\end{equation}
\end{lem}

\begin{figure}
\medskip
\includegraphics [width=12cm,clip]{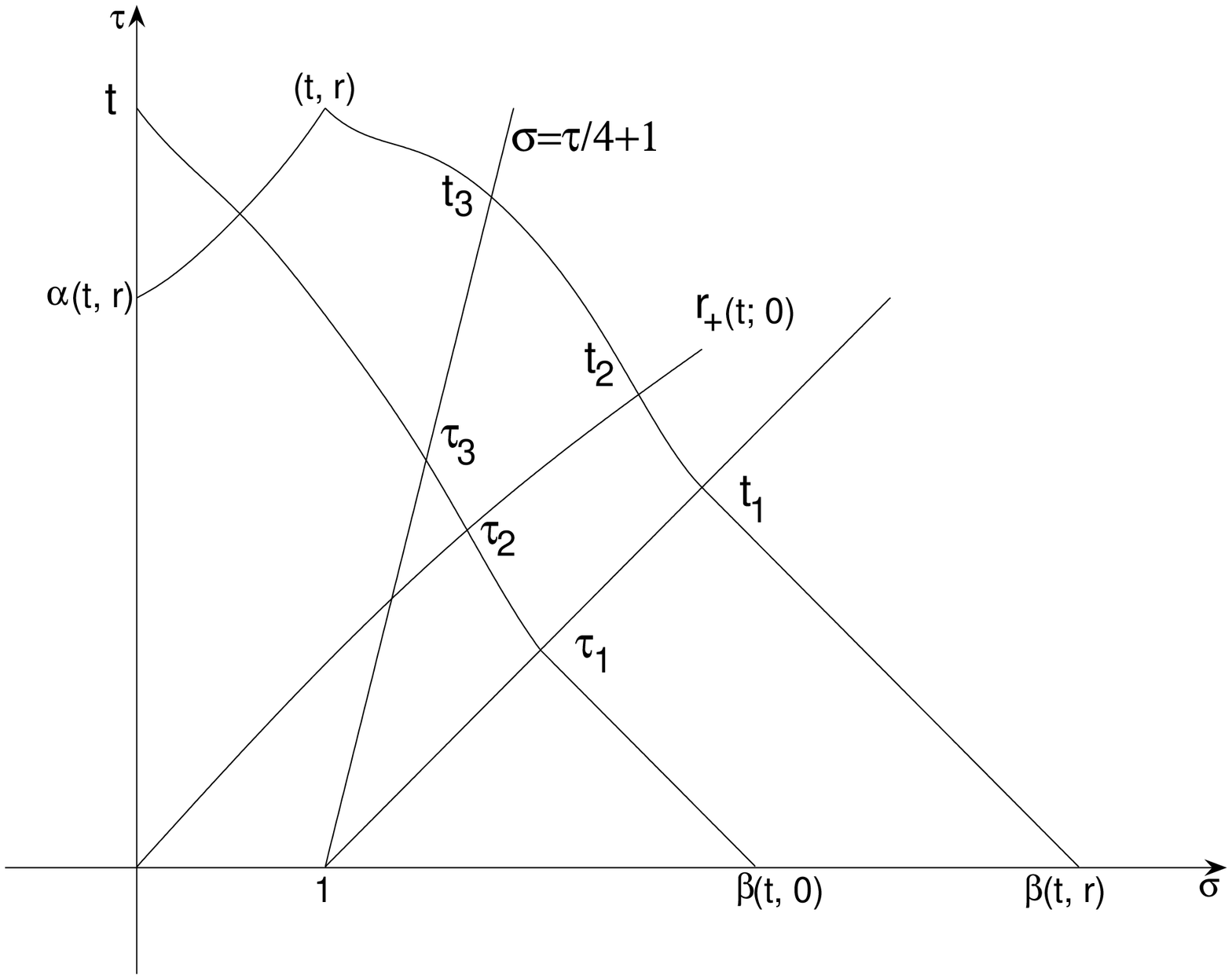}\\
\caption{}
\medskip
\end{figure}

\begin{proof}
First of all, by using the the third line of \eqref{c2} and the
fifth equality of \eqref{c1}, we have
\begin{eqnarray}\nonumber
&&\Big|\frac{d\beta\big(\tau, r_+(\tau; \alpha)\big)}
  {d\tau}\Big|\\\nonumber
&&= \big|\big(\beta_t + a(u)\beta_r\big)\big(\tau, r_+(\tau;
  \alpha)\big)\big|\\\nonumber
&&= \frac{2a(u)\big(\tau, r_+(\tau; \alpha)\big)}{r_{-,
  \beta}\big[\tau; \beta\big(\tau, r_+(\tau;
  \alpha)\big)\big]}\\\nonumber
&&= 2a(u)\big(\tau, r_+(\tau; \alpha)\big)
  \exp\Big\{\int_0^\tau a^\prime(u)u_r
  \Big(s, r_-\big[s; \beta\big(\tau, r_+(\tau; \alpha)\big)
  \big]\Big)ds\Big\}.
\end{eqnarray}
Obviously, the above calculations are also valid with $\alpha$
being replaced by $\gamma$. Thus, we can compute
\begin{eqnarray}\nonumber
&&\int_\alpha^t\frac{1}{\big[1 + \beta\big(\tau,
  r_+(\tau; \alpha)\big)\big]^{2 - \frac{10}{K}}}d\tau\\\nonumber
&&\leq \int_{\beta(\alpha, 0)}^\infty\frac{1}{(1 +
  \beta)^{2 - \frac{10}{K}}}\frac{d\beta}
  {\Big|\frac{d\beta\big(\tau, r_+(\tau;
  \alpha)\big)}{d\tau}\Big|}\\\nonumber
&&\leq \int_{\beta(\alpha, 0)}^\infty\frac{4\exp\Big\{
  C\int^{\tau(\alpha, \beta)}_0\|u_r(s, \cdot)\|_{L^\infty}
  ds\Big\}}{1 + \beta^{2 - \frac{10}{K}}}d\beta,
\end{eqnarray}
and similarly,
\begin{eqnarray}\nonumber
&&\int_0^t\frac{1}{\big[1 + \beta\big(\tau,
  r_+(\tau; \gamma)\big)\big]^{2 - \frac{10}{K}}}d\tau\\\nonumber
&&\leq \int_\gamma^\infty\frac{4\exp\Big\{
  C\int^{\tau(\gamma; \beta)}_0\|u_r(s, \cdot)\|_{L^\infty}
  ds\Big\}}{1 + \beta^{2 - \frac{10}{K}}}d\beta.
\end{eqnarray}
Recall Remark \ref{rem21} and note the first inequality of
\eqref{c6}, we have
\begin{eqnarray}\label{c10}
&&\exp\Big\{C\int^{\tau(\alpha, \beta)}_0\|u_r(s,
  \cdot)\|_{L^\infty}ds\\\nonumber
&&\leq \big(1 + \tau(\alpha; \beta)\big)^{CK^2\epsilon}\\\nonumber
&&\leq 3(1 + \beta)^{CK^2\epsilon},
\end{eqnarray}
Consequently,
\begin{eqnarray}\nonumber
&&\int_\alpha^t\frac{1}{\big[1 + \beta\big(\tau,
  r_+(\tau; \alpha)\big)\big]^{2 - \frac{10}{K}}}d\tau\\\nonumber
&&\leq \int_{\beta(\alpha, 0)}^\infty\frac{12}{(1 +
  \beta)^{2 - \frac{10}{K}- CK^2\epsilon}} d\beta \leq  C.
\end{eqnarray}
The second inequality in \eqref{c8} follows by similar arguments.

To prove the last two inequalities of \eqref{c8}, when $t < t_2$,
we deduce that
\begin{eqnarray}\label{c11}
&&t - t_1 = \int_{t_1}^td\tau\\\nonumber &&\leq \int_{t_1}^t
  \big[a(u)\big(\tau, r_+(\tau; 1)\big) + a(u)\big(\tau, r_-(\tau;
  \beta)\big) \big]d\tau\\\nonumber
&&= \int^t_{t_1}\frac{dr_+(\tau;
  1)}{d\tau}d\tau - \int^t_{t_1}\frac{dr_-(\tau;
  \beta)}{d\tau}d\tau\\\nonumber
&&= r_+(t; 1) - r_+(t_1; 1) - r_-(t; \beta) + r_-(t_1;
  \beta)\\\nonumber
&&= r_+(t; 1) - r_+(t; \gamma) = \int_\gamma^1r_{+, \gamma}(t;
  \mu)d\mu\\\nonumber
&&\leq (1 + t)^{CK^2\epsilon},
\end{eqnarray}
where in the third equality we used \eqref{c7}, and in the last
inequality we used the third equality of \eqref{c1} and Remark
\ref{rem21}. Thus, note that $0 \leq \gamma \leq 1$, we have
\begin{eqnarray}\nonumber
\int_{t_1}^t\frac{1}{\big[1 + \gamma\big(\tau, r_-(\tau;
  \beta)\big)\big]^{2 - \frac{10}{K}}}d\tau
\leq t - t_1 \leq (1 + t)^{CK^2\epsilon},.
\end{eqnarray}

Similarly, when $t \geq t_2$, the combination of the first line of
\eqref{c2} and the first equality of \eqref{c1} gives
\begin{eqnarray}\nonumber
&&\int_{t_2}^t
  \frac{1}{\big[1 + \alpha\big(\tau, r_-(\tau;
  \beta)\big)\big]^{2 - \frac{10}{K}}}d\tau\\\nonumber
&&\leq \int^{\alpha\big(t, r_+(t, \beta)\big)}_0\frac{1}{(1 +
  \alpha)^{2 - \frac{10}{K}}}\frac{1}{\Big|\frac{d\alpha
  \big(\tau, r_-(\tau; \beta)\big)}{d\tau}\Big|}d\alpha\\\nonumber
&&\leq C\int^t_0\frac{C\big(1 + \tau(\alpha,
  \beta)\big)^{CK^2\epsilon}}{(1 +
  \alpha)^{2 - \frac{10}{K}}}d\alpha\\\nonumber
&&\leq C(1 + t)^{CK^2\epsilon}.
\end{eqnarray}
The proof of the Lemma \ref{lem35} is completed.
\end{proof}

\begin{lem}\label{lem36}
Suppose that \eqref{b17} and the assumptions $(H_1)-(H_2)$ in
Theorem \ref{thm11} are satisfied. Then
\begin{equation}\label{c12}
\frac{r_-\big[\tau; \beta(t, r)\big] - r_-\big[\tau;
\beta\big(\alpha(t, r), 0\big)\big]}{r} \leq C(1 +
t)^{CK^2\epsilon},
\end{equation}
holds for any $(t, r)$ with $t \geq t_+(r; 0)$.
\end{lem}
\begin{proof}
A straightforward computation shows that
\begin{eqnarray}\nonumber
&&\Big|\frac{r_-\big[\tau; \beta(t, r)\big] - r_-\big[\tau;
  \beta\big(\alpha(t, r), 0\big)\big]}{r}\Big|\\\nonumber
&&\leq \Big|\frac{\beta(t, r) - \beta\big(\alpha(t, r),
  0\big)}{r}\Big|\sup_{\beta\big(\alpha(t, r), 0\big) \leq \beta
  \leq \beta(t, r)}|r_{-, \beta}(\tau; \beta)|\\\nonumber
&&\leq \frac{1}{r}\int_{\alpha(t, r)}^t\Big|\frac{d\beta\big[\tau,
  r_+\big(\tau; \alpha(t, r)\big)\big]}{d\tau}\Big|d\tau\sup_{
  \beta\big(\alpha(t, r), 0\big) \leq \beta
  \leq \beta(t, r)}|r_{-, \beta}(\tau; \beta)|\\\nonumber
&&\leq \frac{1}{r}\int_{\alpha(t, r)}^t\big|\beta_t +
  a(u)\beta_r\big|\big[\tau, r_+\big(\tau; \alpha(t, r)\big)\big]
  d\tau\sup_{\beta\big(\alpha(t, r), 0\big) \leq \beta
  \leq \beta(t, r)}|r_{-, \beta}(\tau; \beta)|.
\end{eqnarray}
Then \eqref{c12} follows from similar arguments as in proving
Lemma \ref{lem35}.
\end{proof}

\begin{lem}\label{lem37}
For fixed $(t, r)$ with $t \geq t_+(r; 0)$ and $r \leq \frac{t}{4}
+ 1$, let
\begin{equation}\label{c14}
m(\tau, t) = r_-\big(\tau; \beta(t, 0)\big)\quad {\rm for}\quad
\tau \leq \alpha(t, r).
\end{equation}
Then, we have
\begin{equation}\label{c15}
\frac{1}{2} \leq \frac{dm(\tau, t)}{dt} \leq 2 .
\end{equation}
Especially, when $\tau = 0$, we have
\begin{equation}\label{c16}
\frac{1}{2} \leq \frac{d\beta(t, 0)}{dt} \leq 2.
\end{equation}
\end{lem}

\begin{proof}
First of all, we can compute
\begin{eqnarray}\label{c17}
\frac{dm(\tau, t)}{dt} =  r_{-, \beta}\big(\tau; \beta(t,
0)\big)\beta_t(t,  0).
\end{eqnarray}
By the fifth inequality of \eqref{c1} and the last line of
\eqref{c2}, we have
\begin{eqnarray}\label{c18}
&&r_{-, \beta}\big(\tau; \beta(t, 0)\big)\beta_t(t,
  0)\\\nonumber
&&= a\big(u(t, 0)\big)\exp\Big\{- \int_0^\tau
  a^\prime(u)u_r\Big(s,
  r_-\big(s; \beta(t, 0)\big)\Big)ds\\\nonumber
&&\quad +\  \int_0^ta^\prime(u)u_r\Big(s, r_-\big(s; \beta(t,
  0)\big)\Big)ds\Big\}\\\nonumber
&&= a\big(u(t, 0)\big)\exp\Big\{\int_\tau^ta^\prime(u)u_r\Big(s,
  r_-\big(s; \beta(t, 0)\big)\Big)ds\Big\}.
\end{eqnarray}

Define $\tau_1$ and $\tau_3$ by (see figure 2)
\begin{equation}\label{c19}
\begin{cases}
r_+(\tau_1; 1) = r_-\big(\tau_1; \beta(t, 0)\big),\\
r_-\big(\tau_3; \beta(t, 0)\big) = \frac{\tau_3}{4} + 1.
\end{cases}
\end{equation}
The third line of \eqref{c6} gives
\begin{equation}\label{c20}
\tau_1 + 1 \leq \tau_3 + 1 \leq t + 1 \leq 3(\tau_1 + 1).
\end{equation}
Thus, the combination of \eqref{c20} with \eqref{b20}  and
\eqref{b21} gives
\begin{eqnarray}\label{c21}
&&\Big|\int_\tau^ta^\prime(u)u_r\Big(s,
  r_-\big(s; \beta(t, 0)\big)\Big)ds\Big|\\\nonumber
&&\leq \int_{\tau_1}^{\tau_3}\frac{CK\epsilon}{1 +
  s}ds + \Big(\int_{\tau_3}^t(1 + s)^{- 2 +
  \frac{8}{K}}ds\Big)^{\frac{1}{2}}\\\nonumber
&&\quad\times\Big(\int_{\tau_3}^t\Big[(1 + s)^{1 -
  \frac{4}{K}}\sup_{\sigma \leq \frac{s}{4} + 1}
  |u_r(s, \sigma)|\Big]^2ds\Big)^{\frac{1}{2}}\\\nonumber
&&\leq CK\epsilon\ln\frac{1 + \tau_3}{1 + \tau_1} +
  CK^2\epsilon t^{\frac{1}{2}}(1 + \tau_3)^{
  - 1 + \frac{4}{K}}\\\nonumber
&&\leq CK^2\epsilon.
\end{eqnarray}
Finally, noting \eqref{b16} and \eqref{b18}, and combining
\eqref{c17}, \eqref{c18} and \eqref{c21}, we arrive at
\begin{eqnarray}\nonumber
\frac{1}{2} \leq a\big(u(t, 0)\big)\exp\{- CK^2\epsilon\} \leq
\frac{dm(\tau, t)}{dt} \leq a\big(u(t, 0)\big)\exp\{CK^2\epsilon\}
\leq 2
\end{eqnarray}
provided that \eqref{b17} is satisfied.
\end{proof}

\section{Characteristic Method}

In this section we prove the \textit{a priori} estimates
\eqref{b18}, \eqref{b19} and \eqref{b20} in Theorem \ref{thm21} by
using the characteristic method. The strategy of proving the
series of estimates is that by assuming their validity first, we
show these estimates are still valid when the constants $C_j$ are
replaced by $\frac{1}{2}C_j$ for $j = 1, 2, \cdots, 6$.

By using the weighted differential operators $L_\pm$ defined in
\eqref{b12}, we rewrite the nonlinear wave equation in \eqref{a1}
as (see \cite{Lindblad1})
\begin{equation}\label{b13}
\begin{cases}
L_+L_-v = - h(u)L_+uL_+v,\\[-5mm]\\
L_-L_+v = - h(u)L_-uL_-v,
\end{cases}
\end{equation}
where
\begin{equation}\label{b14}
v = ru,
\end{equation}
and
\begin{equation}\label{b15}
h(u) = \frac{a^\prime(u)}{2a(u)}.
\end{equation}
We divide the proof into three steps.

\smallskip
\emph{Step 1. Proof of the a priori estimates for $u$ in
\eqref{b18}:}
\smallskip

First of all, when $t \geq t_+(r; 0)$, noting that $K$ is assumed
to be a big enough constant, and by the last inequality of
\eqref{b19}, we have
\begin{eqnarray}\label{d1}
&&\frac{\int_0^r\big|L_-v(t, \sigma)\big|d\sigma}{r}\\\nonumber &&
  \leq \frac{C_4\epsilon}{r}\int_0^r\frac{1}{\big(1 + \alpha(t,
  \sigma)\big)^{\frac{3}{4}}}d\sigma\\\nonumber
&&\leq \frac{C_4\epsilon}{r}\int_{\alpha(t, r)}^t\frac{1}{\big(1
  + \alpha\big)^{\frac{3}{4}}}\frac{1}{\big|
  \frac{d\alpha(t, \sigma)}{d\sigma}
  \big|}d\alpha.
\end{eqnarray}
By the first line of \eqref{c2} and the first equality of
\eqref{c1}, we can bound the last quantities of the above
inequality by
\begin{eqnarray}\nonumber
\frac{C_4\epsilon}{r}\int_{\alpha(t,
  r)}^t\frac{\big|a\big(u(\alpha, 0)\big)\big|
  \exp\Big\{\int_{\alpha}^ta^\prime(u)u_r\big(s,
  r_+(s, \alpha)\big)ds\Big\}}{\big(1 +
  \alpha\big)^{\frac{3}{4}}}d\alpha.
\end{eqnarray}
Remark \ref{rem21} implies that the above quantity can bounded by
\begin{eqnarray}\nonumber
\frac{2C_4\epsilon(1 +
  t)^{CK^2\epsilon}}{r}\int_{\alpha(t,
  r)}^t\frac{1}{\big(1 +
  \alpha\big)^{\frac{3}{4}}}d\alpha.
\end{eqnarray}
Finally, recall the second inequality in \eqref{c5}, we compute
\begin{eqnarray}\label{d2}
&&\frac{\int_0^r\big|L_-v(t, \sigma)\big|d\sigma}{r}\\\nonumber
&&\leq \frac{8C_4\epsilon(1 + t)^{CK^2\epsilon}}{r}\big[(1 +
  t)^{\frac{1}{4}} - \big(1 + \alpha(t,
  r)\big)^{\frac{1}{4}}\big]\\\nonumber
&&\leq \frac{8C_4\epsilon(1 + t)^{CK^2\epsilon}|t -
  \alpha(t, r)|}{r(1 + t)^{\frac{3}{4}}}\\\nonumber
&&\leq \frac{10C_4\epsilon}{(1 + t)^{\frac{3}{5}}}.
\end{eqnarray}

When $t_+(r; 1) \leq t < t_+(r; 0)$, similarly as in \eqref{c11},
we have
\begin{equation}\label{d3}
|r_+(t; 1) - r_+(t; 0)| \leq (1 + t)^{CK^2\epsilon}.
\end{equation}
Moreover, \eqref{c4} gives
\begin{equation}\label{d4}
(1 - C\epsilon)t \leq r_+(t; 0) \leq (1 + C\epsilon)t.
\end{equation}
Thus, for $t \geq 1$, the combination of \eqref{d2}-\eqref{d4} and
the second inequality of \eqref{b19} gives
\begin{eqnarray}\label{d5}
&&\frac{\int_0^r\big|L_-v(t, \sigma)\big|d\sigma}{r}\\\nonumber
&&= \frac{\int_0^{r_+(t; 0)}\big|L_-v(t, \sigma)\big|d\sigma}{r} +
  \frac{\int_{r_+(t; 0)}^r\big|L_-v(t, \sigma)\big|d\sigma}{r}\\\nonumber
&&\leq \frac{10C_4\epsilon}{(1 + t)^{\frac{3}{5}}} +
  \frac{C_4\epsilon|r_+(t; 1) - r_+(t; 0)|}{r_+(t; 0)}\\\nonumber
&&\leq \frac{12C_4\epsilon}{(1 + t)^{\frac{3}{5}}}.
\end{eqnarray}
For $t \leq 1$, it is obvious that the second inequality in
\eqref{b19} implies
\begin{eqnarray}\label{d6}
\frac{\int_0^r\big|L_-v(t, \sigma)\big|d\sigma}{r} \leq
C_4\epsilon \leq \frac{12C_4\epsilon}{(1 + t)^{\frac{3}{5}}}.
\end{eqnarray}
Finally, combining the first inequality of \eqref{b19} with
\eqref{d2}, \eqref{d5} and \eqref{d6}, we arrive at
\begin{eqnarray}\label{d7}
&&|u(t, r)| = \frac{1}{r}\Big|\int_0^rv_r(t,
  \sigma)d\sigma\Big|\\\nonumber
&&\leq \frac{\sqrt{2}}{r}\int_0^r\big|L_-v(t, \sigma)\big|d\sigma
  + \frac{\sqrt{2}}{r}\int_0^r\big|L_+v(t, \sigma)\big|d\sigma\\\nonumber
&&\leq \frac{18C_4\epsilon}{(1 + t)^{\frac{3}{5}}} +
  \frac{C_3\epsilon(1 + t)^{\frac{1}{K}}}{1 +
  \beta(t, 0)}\\\nonumber
&&\leq \frac{18C_4\epsilon}{(1 + t)^{\frac{3}{5}}} +
  \frac{2C_3\epsilon(1 + t)^{\frac{1}{K}}}{1 +
  t}\\\nonumber
&&\leq \frac{C_1\epsilon}{2(1 + t)^{\frac{3}{5}}}
\end{eqnarray}
by letting
\begin{equation}\label{d8}
C_1 = 40(C_3 + C_4),
\end{equation}
where in the third inequality of \eqref{d7} we used the first
inequality in \eqref{c6}.

Similarly, if we change the estimates \eqref{d1} into
\begin{eqnarray}\nonumber
\frac{\int_0^r\big|L_-v(t, \sigma)\big|d\sigma}{r}
  \leq \frac{C_4\epsilon}{r}\int_0^r\frac{1}{\big(1 + \alpha(t,
  \sigma)\big)^{1 - \frac{1}{K}}}d\sigma,
\end{eqnarray}
we will arrive at
\begin{eqnarray}\nonumber
\frac{\int_0^r\big|L_-v(t, \sigma)\big|d\sigma}{r} \leq
\frac{C_2K\epsilon}{4}(1 + t)^{1 - \frac{2}{K}}
\end{eqnarray}
for $t \geq t_+(r; 0)$, provided that \eqref{b17} holds and
\begin{equation}\label{d9}
C_2 \geq 4(C_3 + C_4).
\end{equation}
The second inequality of \eqref{b18} follows by the same procedure
as in \eqref{d5}-\eqref{d7}.

\bigskip
\emph{Step 2. Proof of the a priori estimates for $L_\pm v$ in
\eqref{b19}:}
\bigskip

Integrating the second equation in \eqref{b13} with respect to
$\tau$ from $0$ to $t$ along the minus characteristics starting
from $\big(0, \beta(t, r)\big)$ yields
\begin{eqnarray}\label{d10}
&&L_+v(t, r) = L_+v(0, \beta) -
  \int_{t_1}^th(u)L_-uL_-v\big(\tau, r_-(\tau;
  \beta)\big)d\tau,
\end{eqnarray}
where $t_1$ is given in \eqref{c7} (see Figure 2).

Define $t_3$ by (see Figure 2)
\begin{eqnarray}\nonumber
\frac{t_3}{4} + 1 = r_-\big(t_3; \beta(t, r)\big).
\end{eqnarray}
Recalling the definition of $t_2$ in \eqref{c9} and using
\eqref{b19}, \eqref{b20}, \eqref{b21} and the last inequalities in
Lemma \ref{lem35}, we compute
\begin{eqnarray}\label{d11}
&&\Big|\int_{t_1}^th(u)L_-uL_-v\big(\tau, r_-(\tau;
  \beta)\big)d\tau\Big|\\\nonumber
&&\leq CK\epsilon^2\int_{t_1}^{t_2}\frac{1}{1 + \tau}d\tau +
  CK\epsilon^2\int_{t_2}^{t_3}\frac{1}{\big[1 + \alpha\big(\tau,
  r_-(\tau; \beta)\big)\big]^{1 - \frac{1}{K}}}\\\nonumber
&&\quad \times \frac{1}{(1 + \tau)\big[1 + \alpha\big(\tau,
  r_-(\tau; \beta)\big)\big]^{1 -
  \frac{2}{K}}}d\tau\\\nonumber
&&\quad +\ C\epsilon\Big(\int_{t_3}^t(1 + \tau)^{-2 +
  \frac{8}{K}}\sup_{\sigma \leq \frac{\tau}{4} + 1}
  \big|(1 + \tau)^{1 - \frac{4}{K}}
  L_-u(\tau, \sigma)\big|^2d\tau\Big)^{\frac{1}{2}}\\\nonumber
&&\quad \times\Big(\int_{t_3}^t\frac{1}{\big[1 + \alpha\big(\tau,
  r_-\big(\tau; \beta)\big)\big]^{2 - \frac{2}{k}}}d\tau\Big)^{\frac{1}{2}}\\\nonumber
&&\leq CK\epsilon^2\frac{t_2 - t_1}{1 + t_1} +
  \frac{CK\epsilon^2(1 + t_3)^{CK^2\epsilon}}{1 + t_2}\\\nonumber
&&\quad +\ CK^2\epsilon^2(1 + t_3)^{- 1 + \frac{4}{K}}(1
  + t_3)^{- \frac{1}{2} + \frac{1}{K}},
\end{eqnarray}
where we used the second line of \eqref{c6} to estimate the
quantity in the fifth line of the above equalities. Finally, with
the aid of \eqref{c11}, the first and third lines of \eqref{c6},
the combination of \eqref{d10} and \eqref{d11} gives
\begin{eqnarray}\label{d12}
&&\big|\beta(t, r)L_+v(r, t)\big|\leq \big|\beta(t,
  r)L_+v\big(0, \beta(t, r)\big)\big|\\\nonumber
&&\quad  +\ \beta(t, r)\int_{t_1}^t\big|h(u)L_-uL_-v\big[\tau,
  r_-\big(\tau; \beta(t, r)\big)\big]\big|d\tau\\\nonumber
&&\leq \sup_{\beta \leq 1}\{|\beta L_+v(0, \beta)|\} +
  CK^2\epsilon^2(1 + t)(1 + t)^{- 1 + CK^2\epsilon}\\\nonumber
&&\leq 2(\|\widetilde{f}^\prime\|_{L^\infty} +
  \|\widetilde{g}\|_{L^\infty}) + CK^2\epsilon^2(1 +
  t)^{CK^2\epsilon}\\\nonumber
&&\leq \frac{C_3\epsilon(1 + t)^{\frac{1}{K}}}{2}.
\end{eqnarray}
where
\begin{eqnarray}\label{d13}
\widetilde{f}(r) = rf(r),\quad \widetilde{g}(r) = rg(r),
\end{eqnarray}
provided that \eqref{b17} holds and
\begin{equation}\label{d14}
\|\widetilde{f}^\prime\|_{L^\infty} + \|\widetilde{g}\|_{L^\infty}
\leq \frac{C_3\epsilon}{4}.
\end{equation}
By Hardy inequality, it is obviously that \eqref{a5} implies
\eqref{d14} by the straight forward calculations:
\begin{eqnarray}\label{d171}
&&\|\widetilde{f}^\prime\|_{L^\infty} +
  \|\widetilde{g}\|_{L^\infty}\\\nonumber
&&\leq C\big(\|\widetilde{f}^{\prime\prime}\|_{L^1} +
  \|\widetilde{g}^\prime\|_{L^1}\big)\\\nonumber
&&\leq C\big(\|\widetilde{f}^{\prime\prime}\|_{L^2} +
  \|\widetilde{g}^\prime\|_{L^2}\big)\\\nonumber
&&\leq C\big(\|f^\prime\|_{L^2} + \|rf^{\prime\prime}\|_{L^2} +
  \|g\|_{L^2} + \|rg^\prime\|_{L^2}\big)\\\nonumber
&&\leq C\big(\|\nabla f\|_{H^1(\mathbb{R}^3)} +
 \|g\|_{H^1(\mathbb{R}^3)}\big).
\end{eqnarray}

Next, we prove the second and third inequalities in \eqref{b19}.
When $t < t_+(r; 0)$, integrating the first equation in
\eqref{b13} along the plus characteristic with respect to $t$
gives
\begin{eqnarray}\nonumber
&&|L_-v(t, r)| = |L_-v(0, \gamma) - \int_0^th(u)L_+uL_+v\big(\tau,
  r_+(\tau; \gamma)\big) d\tau|\\\nonumber
&&\leq 2(\|\widetilde{f}^\prime\|_{L^\infty} +
  \|\widetilde{g}\|_{L^\infty})
  + C\Big(\int_0^t|L_+u|^2\big(\tau, r_+(\tau;
  \gamma)\big)d\tau\Big)^{\frac{1}{2}}\\\nonumber
&&\quad \times \Big(\int_0^t|L_+v|^2\big(\tau, r_+(\tau;
  \gamma)\big)d\tau\Big)^{\frac{1}{2}}.
\end{eqnarray}
By the first inequality in \eqref{b19} and the second inequality
in \eqref{c8}, there holds
\begin{eqnarray}\label{d15}
&&\int_0^t|L_+v\big(\tau, r_+(\tau; \gamma)\big)
  |^2d\tau\\\nonumber
&&\leq C_3^2\epsilon^2\int_0^t\frac{(1 + \tau)^{\frac{2}{K}}}{1
  + \beta^2\big(\tau, r_+(\tau; \gamma)\big)}d\tau\\\nonumber
&&\leq CC_3^2\epsilon^2
\end{eqnarray}
On the other hand, \eqref{b20} and \eqref{b21} gives
\begin{eqnarray}\label{d16}
&&\int_0^t|L_+u\big(\tau, r_+(\tau; \gamma)\big)
  |^2d\tau\\\nonumber
&&\leq \int_0^t\big[(1 + \tau)^{1 - \frac{4}{K}}\sup_{\sigma
  \leq \frac{\tau}{4} + 1}|L_+u(\tau, \sigma)|\big]^2d\tau\\\nonumber
&&\quad +\ \int_0^t\sup_{\sigma >
  \frac{\tau}{4} + 1}|L_+u(\tau, \sigma)|^2d\tau\\\nonumber
&&\leq CK^4\epsilon^2.
\end{eqnarray}
Thus, for $t < t_+(r; 0)$, we have
\begin{eqnarray}\label{d17}
&&|L_-v(t, r)|  \leq \frac{1}{4}C_4\epsilon
\end{eqnarray}
provided that
\begin{equation}\label{d18}
\|\widetilde{f}^\prime\|_{L^\infty} + \|\widetilde{g}\|_{L^\infty}
\leq \frac{C_4\epsilon}{16}.
\end{equation}

Similarly, when $t \geq t_+(r; 0)$, we have
\begin{eqnarray}\label{d19}
&&\alpha|L_-v(t, r)| = \big|\alpha L_-v(\alpha, 0) -
  \alpha\int_\alpha^th(u)L_+uL_+v\big(\tau,
  r_+(\tau; \alpha)\big) d\tau\big|\\\nonumber
&&= \big|\alpha L_+v(\alpha, 0) + \alpha\int_\alpha^t
  h(u)L_+uL_+v\big(\tau, r_+(\tau; \alpha)\big) d\tau\big|\\\nonumber
&&= \big|\alpha L_+v\big(0, \beta(\alpha, 0)\big) - \alpha
  \int_0^\alpha h(u)L_-uL_-v\big[\tau, r_-(\tau; \beta(\alpha, 0)\big)\big]
  d\tau\\\nonumber
&&\quad +\ \alpha\int_\alpha^t h(u)L_+uL_+v\big(\tau, r_+(\tau;
  \alpha)\big) d\tau\big|.
\end{eqnarray}

Define $t_4$ by
\begin{eqnarray}\label{d20}
r_+(t_4; \alpha) = \frac{t_4}{4} + 1.
\end{eqnarray}
By the first inequality in \eqref{b19}, the second inequality in
\eqref{b20} and \eqref{b21}, we have
\begin{eqnarray}\label{d21}
&&\Big|\alpha\int_\alpha^t h(u)L_+uL_+v\big(\tau, r_+(\tau;
  \alpha)\big) d\tau\Big|\\\nonumber
&&\leq C\epsilon\Big(\int_\alpha^{t_4}\big[(1 + \tau)^{1 -
  \frac{4}{K}}\big|L_+u\big(\tau, r_+(\tau; \alpha)\big)\big]^2
  d\tau\Big)^{\frac{1}{2}}\\\nonumber
  &&\quad \times\Big(\int_\alpha^{t_4}\frac{(1 + \tau)^{
  \frac{8}{K}}(1 + \tau)^{\frac{2}{K}}}{\big[1 + \beta\big(\tau,
  r_+(\tau; \alpha)\big)\big]^2}d\tau\Big)^{\frac{1}{2}}\\\nonumber
&&\quad +\ CK\epsilon^2\int^t_{t_4} \frac{\alpha}{(1 + \tau)(1 +
  \alpha)^{1 - \frac{2}{K}}}\frac{(1 + \tau)^{\frac{1}{K}}}{1
  + \beta\big(\tau, r_+(\tau; \alpha)\big)}\big(\tau,
  r_+(\tau; \alpha)\big) d\tau\\\nonumber
&&\leq CK^2\epsilon^2 + CK\epsilon^2\int^t_{t_4}(1 + \tau)^{- 2 +
\frac{3}{K}}d\tau\\\nonumber&&\leq CK^2\epsilon^2.
\end{eqnarray}
On the other hand, \eqref{d11} gives
\begin{eqnarray}\label{d22}
&&\Big|\alpha\int_0^\alpha h(u)L_-uL_-v\big[\tau, r_-(\tau;
  \beta(\alpha, 0)\big)\big]d\tau\Big|\\\nonumber
&&\leq CK^2\epsilon^2\alpha(1 + \alpha)^{- 1 + CK^2\epsilon}.
\end{eqnarray}

The combination of \eqref{d19}-\eqref{d21} gives the third
inequality of \eqref{b19} for $\alpha \geq 1$ provided that
\eqref{d18} holds. Obviously, for $\alpha \leq 1$, the above
procedure is also valid by replacing $\alpha$ by 1. We in fact
have proved \eqref{b19}.

\bigskip
\emph{Step 3. Proof of the a priori estimates for $\partial u$ in
\eqref{b20}:}
\bigskip

By the second inequality of \eqref{b18} and the first and third
inequalities of \eqref{b19}, a straightforward computation shows
\begin{eqnarray}\nonumber
&&\big(1 + \alpha(t, r)\big)^{1 - \frac{2}{K}}
  |\partial u(t, r)|\\\nonumber
&&\leq \big(1 + \alpha(t, r)\big)^{1 -
  \frac{2}{K}}\big(|L_+u(t, r)| + |L_-u(t, r)|\big)\\\nonumber
&&= \big(1 + \alpha(t, r)\big)^{1 -
  \frac{2}{K}}\frac{|rL_+u(t, r)| + |rL_-u(t, r)|}{r}\\\nonumber
&&\leq 4\big(1 + \alpha(t, r)\big)^{1 -
  \frac{2}{K}}\frac{|L_+v(t, r)| + |L_-v(t, r)| + 4|u(t, r)|}{t + 1}\\\nonumber
&&\leq 4\big(1 + \alpha(t, r)\big)^{1 -
  \frac{2}{K}}\frac{1}{t + 1}\Big\{\frac{C_3\epsilon}
  {(1 + t)^{1 - \frac{1}{K}}}\\\nonumber
&&\quad +\ \frac{C_4\epsilon}{\big(1 + \alpha(t,
  r)\big)^{1 - \frac{1}{K}}} + \frac{4C_2K\epsilon}{(1 + t)^{
  1 - \frac{2}{K}}}\Big\}\\\nonumber
&&\leq \frac{C_5K\epsilon}{2(1 + t)}
\end{eqnarray}
provided that $C_5 \geq 32C_2$ and $K$ is big enough, since $r
> \frac{t}{4} + 1$ and $t \geq t_+(r, 0)$. This gives the second
inequality of \eqref{b20}.

Similarly, we prove the first inequality of \eqref{b20} as
follows:
\begin{eqnarray}\nonumber
&&|\partial u(t, r)|\\\nonumber &&\leq 4\frac{|L_+v(t, r)| +
  |L_-v(t, r)| + 4|u(t, r)|}{t + 1}\\\nonumber
&&\leq \frac{4}{t + 1}\Big\{\frac{C_3\epsilon}
  {(1 + t)^{1 - \frac{1}{K}}}\\\nonumber
&&\quad +\ C_4\epsilon + \frac{4C_2K\epsilon}{(1 + t)^{
  1 - \frac{2}{K}}}\Big\}\\\nonumber
&&\leq \frac{C_5K\epsilon}{2(1 + t)}
\end{eqnarray}
provided that $C_5 \geq 32C_2$.

\section{Generalized Energy Method}

To prove the weighted estimate in \eqref{b21}, we need more decay
inside the cone $\mathcal{D}$. For this purpose, we invoke
Klainerman's vector fields and the generalized energy method.

The so-called angular momentum operators are the vector fields
\begin{equation}\nonumber
\Omega = (\Omega_1, \Omega_2, \Omega_3) = x \wedge \nabla,
\end{equation}
where $\wedge$ is the usual vector cross product. The scaling
operator $S$ is defined by
\begin{equation}\nonumber
S = t\partial_t + r\partial_r.
\end{equation}
The Lorentz operators are the vector fields
\begin{equation}\nonumber
\Gamma_j = x_j\partial_t + t\partial_j
\end{equation}
for $j = 1, 2, 3$.

We further introduce the notations
\begin{equation}\label{e1}
\begin{cases}
\widetilde{\Gamma}_1 = (t + r)(\partial_t + \partial_r),\\
\widetilde{\Gamma}_2 = (t - r)(\partial_t - \partial_r),
\end{cases}
\end{equation}
and
\begin{equation}\nonumber
\Gamma_0 = t\partial_r + r\partial_t.
\end{equation}
It is easy to see that
\begin{equation}\label{e2}
\begin{cases}
\widetilde{\Gamma}_1 = S + \Gamma_0,\\
\widetilde{\Gamma}_2 = S - \Gamma_0.
\end{cases}
\end{equation}

For radially symmetric functions $u(t, x) = u(t, |x|)$ defined on
$\mathbb{R} \times \mathbb{R}^3$, there holds
\begin{equation}\nonumber
\Gamma_ju = \frac{x_j}{r}\Gamma_0u,\quad \Omega u = 0.
\end{equation}
Hence, for $\Gamma$ being any of the vector fields
\begin{equation}\nonumber
S,\quad \Gamma_j, \quad \Omega_{ij},\quad i, j = 1, 2, 3,
\end{equation}
there holds
\begin{equation}\label{e3}
\|\Gamma u(t, \cdot) \|_{L^\infty(\mathbb{R}^3)} \leq
\frac{1}{2}\big(\|\widetilde{\Gamma}_1 u(t, \cdot)\|_{L^\infty} +
\|\widetilde{\Gamma}_2 u(t, \cdot)\|_{L^\infty}\big).
\end{equation}

Applying $\Gamma$ to the quasi-linear wave equation in \eqref{a1},
we have
\begin{eqnarray}\nonumber
&&(\Gamma u)_{tt} - a^2(u)\Delta(\Gamma u)\\\nonumber && =
  [\partial_t^2, \Gamma]u + \Gamma\partial_t^2u - a^2(u)
  \big([\Delta, \Gamma]u + \Gamma\Delta u\big)\\\nonumber
&&= [\partial_t^2, \Gamma]u - a^2(u)[\Delta, \Gamma]u
  + 2a(u)a^\prime(u)\Gamma u\Delta u\\\nonumber
&&= 2a(u)a^\prime(u)\Gamma u\Delta u\\\nonumber &&\quad +\
\begin{cases}
0,\quad {\rm for}\quad \Gamma =
S\ {\rm or}\ \Omega_{ij},\\
2\big(1 - a^2(u)\big)u_{tx_j},\quad {\rm for}\quad \Gamma =
\Gamma_j,
\end{cases}
\quad i, j = 1, 2, 3.
\end{eqnarray}
The standard energy estimate (see \eqref{b2} for the definition of
$E_s$) gives
\begin{eqnarray}\label{e4}
&&\frac{d}{dt}E_1\big[\Gamma u(t)\big]) \leq C\|\partial u(t,
  \cdot)\|_{L^\infty(\mathbb{R}^3)}E_1\big[\Gamma
  u(t)\big]\\\nonumber
&&+\ C\Big(\|\Delta u\Gamma u(t, \cdot)\|_{L^2(\mathbb{R}^3)} +
  \|u(t, \cdot)\|_{L^\infty(\mathbb{R}^3)}E_2\big[u(t)
  \big]^{\frac{1}{2}}\Big)E_1\big[\Gamma u(t)\big]^{\frac{1}{2}}.
\end{eqnarray}

Noting \eqref{e1}, \eqref{e2} and \eqref{e3}, we compute that
\begin{eqnarray}\label{e5}
&&\|\Delta u\Gamma u(t, \cdot)\|_{L^2(\mathbb{R}^3)}\\\nonumber
&&\leq \|\Delta u\Gamma u(t, x)\|_{L^2(|x| > \frac{t}{4}
  + 1)} + \|\Delta u\Gamma u(t, x)\|_{L^2(|x| \leq
  \frac{t}{4} + 1)}\\\nonumber
&&\leq \|\Delta u(t, x)\|_{L^2(|x| > \frac{t}{4} +
1)}\|\widetilde{\Gamma}
  u(t, r)\|_{L^\infty(r > \frac{t}{4} + 1)}\\\nonumber
&&\quad +\ \big\|\big(1 + \big|t - |x|\big|\big)\Delta u(t,
x)\big\|_{L^2(|x| \leq \frac{t}{4} +
  1)}\Big\|\frac{\widetilde{\Gamma} u(t, r)}{1 + |t - r|}\Big\|_{
  L^\infty(r \leq \frac{t}{4} + 1)}\\\nonumber
&&\leq E_2\big[u(t)\big]^{\frac{1}{2}}\|\widetilde{\Gamma} u(t, r)
  \|_{L^\infty(r > \frac{t}{4} + 1)}\\\nonumber
&&\quad +\  C\Big(E_2\big[u(t)\big]^{\frac{1}{2}} +
  E_1\big[\Gamma u(t)\big]^{\frac{1}{2}}\Big)\|\partial
  u(t, r)\|_{L^\infty(r \leq \frac{t}{4} + 1)}.
\end{eqnarray}
By the second and third inequalities of \eqref{c5}, a straight
forward computation yields
\begin{eqnarray}\nonumber
&&\|\widetilde{\Gamma}_1u(t, \cdot)\|_{L^\infty} = \|(t +
  r)(\partial_t + \partial_r)u(t, \cdot)\|_{L^\infty}\\\nonumber &&\leq
  \|(t - r)(\partial_t  + \partial_r)u(t, r)\|_{L^\infty}
  + 2\|r(\partial_t  + \partial_r) u(t, r)\|_{L^\infty}\\\nonumber
&&\leq \|\alpha(t, r)\partial u(t, r)\|_{L^\infty\big(r \leq
  r_+(t; 0)\big)} + \|\gamma(t, r)\partial u(t, r)\|_{L^\infty\big(r
  > r_+(t; 0)\big)}\\\nonumber
&&\quad +\ (2 + C\epsilon)\|r(\partial_t + \partial_r)
  u(t, r)\|_{L^\infty}\\\nonumber
&&\leq \|\alpha(t, r)\partial u(t, r)\|_{L^\infty\big(r \leq
  r_+(t; 0)\big)} + \|\partial u(t, r)\|_{L^\infty\big(r
  > r_+(t; 0)\big)}\\\nonumber
&&\quad +\ 3\|(\partial_t + \partial_r) v(t, \cdot)\|_{L^\infty} +
  3\|u(t, \cdot)\|_{L^\infty}\\\nonumber
&&\leq \|\alpha(t, r)\partial u(t, r)\|_{L^\infty\big(r \leq
  r_+(t; 0)\big)} + \|\partial u(t, r)\|_{L^\infty
  \big(r > r_+(t; 0)\big)}\\\nonumber
&&\quad +\ 3\|L_+v(t, \cdot)\|_{L^\infty} + C\|u\partial v(t,
  \cdot)\|_{L^\infty} + 3\|u(t, \cdot)\|_{L^\infty}.
\end{eqnarray}
Recall \eqref{b20}, \eqref{b19} and \eqref{b18}, we have
\begin{eqnarray}\label{e6}
&&\|\alpha(t, r)\partial u(t, r)\|_{L^\infty\big(\frac{t}{4} + 1
  \leq r \leq r_+(t; 0)\big)} + \|\partial u(t, r)\|_{L^\infty
  \big(r > \max\{r_+(t; 0), \frac{t}{4} + 1\}\big)}\\\nonumber
&&\quad +\ 3\|L_+v(t, \cdot)\|_{L^\infty} + C\|u\partial v(t,
  \cdot)\|_{L^\infty} + 3\|u(t, \cdot)\|_{L^\infty}\\\nonumber
&&\leq \frac{CK\epsilon\big(1 + \alpha(t, r)\big)^{
  \frac{2}{K}}}{1 + t} + \frac{CK\epsilon}{1 + t} + C\epsilon(1 + t)^{
  - 1 + \frac{2}{K}}\\\nonumber
&&\quad +\ CK\epsilon^2(1 + t)^{- 1 + \frac{2}{K}} +
  CK\epsilon(1 + t)^{- 1 + \frac{2}{K}}\\\nonumber
&&\leq CK\epsilon(1 + t)^{- 1 + \frac{2}{K}},
\end{eqnarray}
which means that
\begin{eqnarray}\label{e7}
\|\widetilde{\Gamma}_1u(t, \cdot)\|_{L^\infty(r > \frac{t}{4} +
1)} \leq CK\epsilon(1 + t)^{- 1 + \frac{2}{K}}.
\end{eqnarray}
Similarly, by the second inequality in \eqref{c5}, one has
\begin{eqnarray}\label{e8}
&&\|\widetilde{\Gamma}_2u(t, \cdot)\|_{L^\infty(r > \frac{t}{4}
  + 1)}\\\nonumber
&&= \|(t - r)(\partial_t - \partial_r)u(t, \cdot)
  \|_{L^\infty(r > \frac{t}{4} + 1)}\\\nonumber
&&\leq \|\alpha(t, r)\partial u(t, r)\|_{L^\infty\big(\frac{t}{4}
  + 1 < r \leq r_+(t; 0)\big)}\\\nonumber
&&\quad +\ CK^2\epsilon(1 + t)^{\frac{2}{K}}
  \|\partial u(t, \cdot)\|_{L^\infty(r > \frac{t}{4} + 1)}\\\nonumber
&&\quad +\ \big\|\gamma(t, r)\partial u(t, \cdot)\big\|_{L^\infty
  \big(r > \max\big\{r_+(t; 0), \frac{t}{4} + 1\big\}\big)}\\\nonumber
&&\leq CK\epsilon(1 + t)^{- 1 + \frac{2}{K}}.
\end{eqnarray}
Inserting \eqref{e7} and \eqref{e8} into \eqref{e5} gives
\begin{eqnarray}\label{e9}
&&\|\Delta u\Gamma u(t, \cdot)\|_{L^2(\mathbb{R}^3)}\\\nonumber
&&\leq C\big[K\epsilon(1 + t)^{- 1 + \frac{2}{K}}
  + \|\partial u(t, r)\|_{L^\infty(r \leq \frac{t}{4} + 1)}\big]
  E_2\big[u(t)\big]^{\frac{1}{2}}\\\nonumber &&\quad +\
  CE_1\big[\Gamma u(t)\big]^{\frac{1}{2}}\|\partial
  u(t, r)\|_{L^\infty(r \leq \frac{t}{4} + 1)}
\end{eqnarray}

Now we use \eqref{e9}, the second inequality in \eqref{b18} and
the basic energy estimate in Remark \ref{rem21} to improve
\eqref{e4} as
\begin{eqnarray}\nonumber
&&\frac{d}{dt}E_1\big[\Gamma u(t)\big]^{\frac{1}{2}} \leq
  C\|\partial u(t, \cdot)\|_{L^\infty(\mathbb{R}^3)} E_1\big[\Gamma
  u(t)\big]^{\frac{1}{2}}\\\nonumber
&&\quad +\ CE_2\big[u(0)\big]^{\frac{1}{2}}(1 + t)^{CK^2\epsilon}
  \big[K\epsilon (1 + t)^{- 1 + \frac{2}{K}} +
  \|\partial u(t, r)\|_{L^\infty(r \leq \frac{t}{4} + 1)}\big].
\end{eqnarray}
Gronwall's inequality gives
\begin{eqnarray}\label{e10}
&&E_1\big[\Gamma u(t)\big]^{\frac{1}{2}}\\\nonumber &&\leq
  \exp\Big\{\int_0^tC\|\partial u(\tau, \cdot)\|_{L^\infty
  (\mathbb{R}^3)}d\tau\Big\}\Big[CE_1\big(\Gamma
  u(0)\big)^{\frac{1}{2}} + CE_2\big[u(0)\big]^{\frac{1}{2}}\\\nonumber
&&\quad \times\ \int_0^t(1 + s)^{CK^2\epsilon}\big[K\epsilon(1 +
  s)^{- 1 + \frac{2}{K}} + \|\partial u(s, \rho)\|_{L^\infty(\rho
  \leq \frac{s}{4} + 1)}\big]ds\Big]\\\nonumber
&&\leq CE_2\big[u(0)\big]^{\frac{1}{2}}(1 + t)^{
  CK^2\epsilon}\Big\{K\epsilon(1 + t)^{\frac{5}{2K}}
  + \Big(\int_0^t(1 + s)^{- 2 + \frac{8}{K} + CK^2\epsilon}\Big)^{\frac{1}{2}}\\\nonumber
&&\quad \times\ \Big(\int_0^t\big[(1 + s)^{1 -
  \frac{4}{K}}\|\partial
  u(s, \rho)\|_{L^\infty(\rho \leq \frac{s}{4} +
  1)}\big]^2ds\Big)^{\frac{1}{2}}\Big\}\\\nonumber
&&\leq CE_2\big[u(0)\big]^{\frac{1}{2}}(1 + t)^{
  \frac{3}{K}}.
\end{eqnarray}
The combination of \eqref{e1}, \eqref{e2} and \eqref{e10} gives
\begin{eqnarray}\label{e11}
&&\|\partial^2u\|_{L^2\big(\mathbb{R}^3(|x| \leq \frac{t}{4} +
1)\big)} \leq \frac{C}{t +
  1}\|\Gamma\partial u\|_{L^2(\mathbb{R}^3)}\\\nonumber
&&\leq \frac{C}{t + 1}\big[E_1\big(\Gamma u(t)
  \big)^{\frac{1}{2}} + \|\partial u\|_{L^2(\mathbb{R}^3)}\big]\\\nonumber
&&\leq \frac{C}{t + 1}\big[E_1\big(\Gamma u(t)
  \big)^{\frac{1}{2}} + E_2\big(u(t)
  \big)^{\frac{1}{2}}\big]\\\nonumber
&&\leq CE_2\big(u(0)\big)^{\frac{1}{2}}(t + 1)^{- 1 +
  \frac{3}{K}}.
\end{eqnarray}

\section{Estimates for Maximal Functions}

This section proves the weighted estimate \eqref{b21}. We will
inevitably encounter the estimation for the second derivatives of
$u$, which will appear when one differentiates equation \eqref{f5}
to get the expression for $\partial u$ (see $A_{11}$, $A_{15}$,
$B_{11}$ and $B_{15}$ in \eqref{f7} and \eqref{f9}) and intends to
use Hardy-Littlewood Maximal functions to eliminate the
singularities of $\partial u$ as $r \rightarrow 0$ (see $A_8$ and
$A_{13}$ in \eqref{f7}). The energy decay property \eqref{e11}
inside the cone $\mathcal{D}$ plays an essential role in
estimating $\partial^2u$.

Another difficulty lies in the inefficiency of decay of $\partial
u$ outside the cone $\mathcal{D}$. To overcome this difficulty, we
rewrite the quasi-linear wave equation as the special form
\eqref{f2}, which takes advantages of the careful estimation for
the accumulation of the quantity $(1 + \alpha)^{- 2 +
\frac{10}{K}}$ along minus characteristics, which combining with
the special form \eqref{f2} provides us extra decay for $\partial
u$ outside the cone $\mathcal{D}$.

As the proof is rather long, we shall divide it into four steps.

\bigskip
\textit{Step 1. Derivation of the expressions for $u_t$ and
$u_r$.}
\bigskip

First of all, by using the weighted operators defined in
\eqref{b27}, we can rewrite the quasi-linear wave equation in
\eqref{a1} as
\begin{equation}\label{f1}
M_-M_+v = - \frac{a^\prime(u)rM_-uu_t}{a^2(u)}.
\end{equation}
Then, we further rewrite the above equation as
\begin{eqnarray}\label{f2}
&&M_-M_+[(t - r) v]\\\nonumber &&= M_-[(t - r) M_+v] +
  M_-\Big[\frac{\big(1 - a(u)\big)v}{a(u)}\Big]\\\nonumber
&&= (t - r)M_-M_+v + \Big[\frac{1 + a(u)}{a(u)}\Big]M_+v +
  M_-\Big[\frac{\big(1 - a(u)\big)v}{a(u)}\Big]\\\nonumber
&&= - \frac{a^\prime(u)(t - r)rM_-uu_t}{a^2(u)} + \Big[\frac{1 +
  a(u)}{a(u)}\Big]M_+v + M_-\Big[\frac{\big(1
  - a(u)\big)v}{a(u)}\Big].
\end{eqnarray}
Integrating the above equation \eqref{f2} with respect to $\tau$
from $0$ to $t$ along the minus characteristic starting from
$\big(0, \beta(t, r)\big)$ yields
\begin{eqnarray}\label{f3}
&&\big(M_+[(t - r) v]\big)(t, r)\\\nonumber &&= \big(M_+[(t - r)
  v]\big)(0, \beta)\\\nonumber
&&\quad +\ \Big[\frac{\big(1
  - a(u)\big)v}{a(u)}\Big](t, r) - \Big[\frac{\big(1
  - a(u)\big)v}{a(u)}\Big](0, \beta)\\\nonumber
&&\quad +\ \int_0^t\big[\big(1 + a(u)\big)M_+v\big]\big(\tau,
  r_-(\tau; \beta)\big)d\tau\\\nonumber
&&\quad -\ \int_0^t\frac{a^\prime(u)(t - r)rM_-uu_t}{a(u)}
  \big(\tau, r_-(\tau; \beta)\big)d\tau.
\end{eqnarray}
For any $(t, r)$ with $t \geq t_+(r; 0)$, $0 \leq \sigma \leq r$,
\eqref{f3} gives
\begin{eqnarray}\label{f4}
&&\frac{d\Big([(t - r) v]\big(t_+(\sigma; \alpha),
  \sigma\big)\Big)}{d\sigma}\\\nonumber
&& = \big(M_+[(t - r)v]\big)\big[0, \beta\big(t_+(\sigma; \alpha),
  \sigma\big)\big]\\\nonumber
&&\quad +\ \Big[\frac{\big(1 - a(u)\big)v}{a(u)}\Big]
  \big(t_+(\sigma; \alpha), \sigma\big) -
  \Big[\frac{\big(1 - a(u)\big)v}{a(u)}\Big]\Big(0,
  \beta\big(t_+(\sigma; \alpha), \sigma\big)\Big)\\\nonumber
&&\quad +\ \int_0^{t_+(\sigma; \alpha)}\big[\big(1 +
  a(u)\big)M_+v\big]\Big(s, r_-\big[s; \beta\big(t_+(\sigma;
  \alpha), \sigma\big)\big]\Big)ds\\\nonumber
&&\quad -\ \int_0^{t_+(\sigma; \alpha)}\frac{a^\prime(u)(t -
  r)rM_-uu_t}{a(u)}\Big(s, r_-\big[s; \beta\big(t_+(\sigma;
  \alpha), \sigma\big)\big]\Big)ds.
\end{eqnarray}
Integrating the above equation with respect to $\sigma$ from $0$
to $r$ gives
\begin{eqnarray}\nonumber
&&[r(t - r)u](t, r)\\\nonumber && = \int_0^r\big(M_+[(t -
  r)v]\big)\big[0, \beta\big(t_+(\sigma; \alpha),
  \sigma\big)\big]d\sigma\\\nonumber
&&\quad +\ \int_0^r\Big[\frac{\big(1 - a(u)\big)v}{a(u)}\Big]
  \big(t_+(\sigma; \alpha), \sigma\big)d\sigma\\\nonumber
&&\quad -\ \int_0^r\Big[\frac{\big(1 - a(u)\big)v}{a(u)}\Big]
  \Big(0, \beta\big(t_+(\sigma; \alpha),
  \sigma\big)\Big)d\sigma\\\nonumber
&&\quad +\ \int_0^r\int_0^{t_+(\sigma; \alpha)}\big[\big(1 +
  a(u)\big)M_+v\big]\Big(s, r_-\big[s; \beta\big(t_+(\sigma;
  \alpha), \sigma\big)\big]\Big)dsd\sigma\\\nonumber
&&\quad -\ \int_0^r\int_0^{t_+(\sigma; \alpha)}\frac{a^\prime(u)(t
  - r)rM_-uu_t}{a(u)}\Big(s, r_-\big[s; \beta\big(t_+(\sigma;
  \alpha), \sigma\big)\big]\Big)dsd\sigma,
\end{eqnarray}
which is equivalent to
\begin{eqnarray}\label{f5}
&&(t - r)u(t, r) = \frac{1}{r}\int_0^r\big(M_+[(t -
  r)v]\big)\big[0, \beta\big(t_+(\sigma; \alpha),
  \sigma\big)\big]d\sigma\\\nonumber
&&\quad +\ \frac{1}{r}\int_0^r\Big[\frac{\big(1 - a(u)\big)v}
  {a(u)}\Big]\big(t_+(\sigma; \alpha), \sigma\big)d\sigma\\\nonumber
&&\quad -\ \frac{1}{r}\int_0^r\Big[\frac{\big(1 - a(u)\big)v}
  {a(u)}\Big]\Big(0, \beta\big(t_+(\sigma; \alpha),
  \sigma\big)\Big)d\sigma\\\nonumber
&&\quad +\ \frac{1}{r}\int_0^r\int_0^{t_+(\sigma;
  \alpha)}\big[\big(1 + a(u)\big)M_+v\big]\Big(s, r_-\big[s;
  \beta\big(t_+(\sigma;
  \alpha), \sigma\big)\big]\Big)dsd\sigma\\\nonumber
&&\quad -\ \frac{1}{r}\int_0^r\int_0^{t_+(\sigma; \alpha)}
  \frac{a^\prime(u)(t - r)rM_-uu_t}{a(u)}\Big(s, r_-\big[s;
  \beta\big(t_+(\sigma; \alpha), \sigma\big)\big]\Big)dsd\sigma,
\end{eqnarray}

Differentiating the above equality \eqref{f5} with respect to $r$
yields
\begin{eqnarray}\nonumber
&&[(t - r)u_r](t, r) - u(t, r) =\\\nonumber &&\quad -\
  \frac{1}{r^2}\int_0^r\big(M_+[(t - r)v]\big)\big[0,
  \beta\big(t_+(\sigma; \alpha),\sigma\big)\big]d\sigma\\\nonumber
&&\quad +\ \frac{\big(M_+[(t - r)v]\big)(0, \beta)}{r}\\\nonumber
&&\quad +\ \frac{1}{r}\int_0^r\big(M_+[(t - r)v]\big)_r\big[0,
  \beta\big(t_+(\sigma; \alpha),\sigma\big)\big]
  \frac{d}{dr}\beta\Big(t_+\big(\sigma; \alpha(t,
  r)\big),\sigma\Big)d\sigma\\\nonumber
&&\quad -\ \frac{1}{r^2}\int_0^r\Big[\frac{\big(1 - a(u)\big)v}
  {a(u)}\Big]\big(t_+(\sigma; \alpha), \sigma\big)d\sigma\\\nonumber
&&\quad +\ \frac{1}{r^2}\int_0^r\Big[\frac{\big(1 - a(u)\big)v}
  {a(u)}\Big]\Big(0, \beta\big(t_+(\sigma; \alpha),
  \sigma\big)\Big)d\sigma\\\nonumber
&&\quad +\ \frac{1}{r}\Big\{\Big[\frac{\big(1 - a(u)\big)v}
  {a(u)}\Big](t, r) - \Big[\frac{\big(1 - a(u)\big)v}
  {a(u)}\Big](0, \beta)\Big\}\\\nonumber
&&\quad +\ \frac{1}{r}\int_0^r\Big[\frac{\big(1 - a(u)\big)v}
  {a(u)}\Big]_t\big(t_+(\sigma; \alpha), \sigma\big)\frac{
  dt_+\big(\sigma; \alpha(t, r)\big)}{dr}d\sigma\\\nonumber
&&\quad -\ \frac{1}{r}\int_0^r\Big[\frac{\big(1 - a(u)\big)v}
  {a(u)}\Big]_r\Big(0, \beta\big(t_+(\sigma; \alpha),
  \sigma\big)\Big)\frac{d\beta\Big(t_+\big(\sigma; \alpha(t, r)\big),
  \sigma\Big)}{dr}d\sigma\\\nonumber
&&\quad -\ \frac{1}{r^2}\int_0^r\int_0^{t_+(\sigma;
  \alpha)}\big[\big(1 + a(u)\big)M_+v\big]
  \Big(s, r_-\big[s; \beta\big(t_+(\sigma;
  \alpha), \sigma\big)\big]\Big)dsd\sigma\\\nonumber
&&\quad +\ \frac{1}{r}\int_0^t\big[\big(1 +
  a(u)\big)M_+v\big]\big(s, r_-(s; \beta)\big)ds\\\nonumber
&&\quad +\ \frac{1}{r}\int_0^r\big[\big(1 + a(u)\big)
  M_+v\big]\big(t_+(\sigma; \alpha), \sigma\big)\frac{d}{dr}
  t_+\big(\sigma; \alpha(t, r)\big)d\sigma\\\nonumber
&&\quad +\ \frac{1}{r}\int_0^r\int_0^{t_+(\sigma;
  \alpha)}\big[\big(1 + a(u)\big)M_+v\big]_r
  \Big(s, r_-\big[s; \beta\big(t_+(\sigma;
  \alpha), \sigma\big)\big]\Big)dsd\sigma\\\nonumber
&&\quad \times\ \frac{d}{dr}r_-\Big[s;
  \beta\Big(t_+\big(\sigma; \alpha(t, r)\big),
  \sigma\Big)\Big]dsd\sigma\\\nonumber
&&\quad + \frac{1}{r^2}\int_0^r\int_0^{t_+(\sigma;
  \alpha)}\frac{a^\prime(u)(t - r) rM_-uu_t}{a(u)}\Big(s,
  r_-\big[s; \beta\big(t_+(\sigma; \alpha),
  \sigma\big)\big]\Big)dsd\sigma\\\nonumber
&&\quad -\ \frac{1}{r}\int_0^t \frac{a^\prime(u)(t - r)
  rM_-uu_t}{a(u)}\big(s, r_-(s; \beta)\big)ds\\\nonumber
&&\quad -\ \frac{1}{r}\int_0^r\frac{a^\prime(u)(t - r)rM_-uu_t}
  {a(u)}\big(t_+(\sigma; \alpha), \sigma\big)\frac{d}{dr}
  t_+\big(\sigma; \alpha(t, r)\big)d\sigma\\\nonumber
&&\quad -\ \frac{1}{r}\int_0^r\int_0^{t_+(\sigma;
  \alpha)}\Big[\frac{a^\prime(u)(t - r)rM_-uu_t}{a(u)}\Big]_r\Big(s,
  r_-\big[s; \beta\big(t_+(\sigma; \alpha),
  \sigma\big)\big]\Big)\\\nonumber
&&\quad \times\ \frac{d}{dr}r_-\Big[s;
  \beta\Big(t_+\big(\sigma; \alpha(t, r)\big), \sigma\Big)\Big]dsd\sigma.
\end{eqnarray}
The above equality can be rewritten as
\begin{eqnarray}\label{f6}
[(t - r) u_r](t, r) = \sum_{i = 1}^{15} A_i,
\end{eqnarray}
where
\begin{equation}\label{f7}
\begin{cases}
A_1 = u(t, r),\\[-2mm]\\
A_2 = \frac{1}{r}\int_0^r\frac{\big(M_+[(t - r)v]\big)(0, \beta) -
  \big(M_+[(t - r)v]\big)\big[0, \beta\big(t_+(\sigma; \alpha),\sigma\big)
  \big]}{r}d\sigma,\\[-2mm]\\
A_3 = \frac{1}{r}\int_0^r\big(M_+[(t - r)v]\big)_r\big[0,
  \beta\big(t_+(\sigma; \alpha),\sigma\big)\big]
  \frac{d}{dr}\beta\Big(t_+\big(\sigma; \alpha(t, r)\big),
  \sigma\Big)d\sigma,\\[-2mm]\\
A_4 = - \frac{1}{r}\int_0^r\frac{1}{r}\Big\{\Big[\frac{\big(1 -
  a(u)\big)v}{a(u)}\Big]\big(t_+(\sigma; \alpha), \sigma\big)
  - \Big[\frac{\big(1 - a(u)\big)v}{a(u)}\Big](t, r)\Big\}d\sigma,\\[-2mm]\\
A_5 = \frac{1}{r}\int_0^r\frac{1}{r}\Big\{\Big[\frac{\big(1 -
  a(u)\big)v}{a(u)}\Big]\Big(0, \beta\big(t_+(\sigma; \alpha),
  \sigma\big)\Big) - \Big[\frac{\big(1 - a(u)\big)v}
  {a(u)}\Big](0, \beta)\Big\}d\sigma,\\[-2mm]\\
A_6 = \frac{1}{r}\int_0^r\Big[\frac{\big(1 - a(u)\big)v}
  {a(u)}\Big]_t\big(t_+(\sigma; \alpha), \sigma\big)\frac{
  dt_+\big(\sigma; \alpha(t, r)\big)}{dr}d\sigma,\\[-2mm]\\
A_7 =  - \frac{1}{r}\int_0^r\Big[\frac{\big(1 - a(u)\big)v}
  {a(u)}\Big]_r\Big(0, \beta\big(t_+(\sigma; \alpha),
  \sigma\big)\Big)\frac{d\beta\Big(t_+\big(\sigma; \alpha(t, r)\big),
  \sigma\Big)}{dr}d\sigma,\\[-2mm]\\
A_8 = \frac{1}{r}\int_0^r\frac{1}{r}
  \Big\{\int_0^{t_+(\sigma; \alpha)}\big[\big(1 + a(u)\big)M_+v\big]
  \big(s, r_-(s; \beta)\big)ds,\\[-2mm]\\
\quad\quad  - \int_0^{t_+(\sigma;
  \alpha)}\big[\big(1 + a(u)\big)M_+v\big]
  \Big(s, r_-\big[s; \beta\big(t_+(\sigma;
  \alpha), \sigma\big)\big]\Big)ds\Big\}d\sigma,\\[-2mm]\\
A_9 = \frac{1}{r}\int_0^r\frac{1}{r}
  \int_{t_+(\sigma; \alpha)}^t\big[\big(1 + a(u)\big)M_+v\big]
  \big(s, r_-(s; \beta)\big)dsd\sigma,\\[-2mm]\\
A_{10} = \frac{1}{r}\int_0^r\big[\big(1 + a(u)\big)M_+v\big]
  \big(t_+(\sigma; \alpha), \sigma\big)\frac{d}{dr}
  t_+\big(\sigma; \alpha(t, r)\big)d\sigma,\\[-2mm]\\
A_{11} = \frac{1}{r}\int_0^r\int_0^{t_+(\sigma;
  \alpha)}\big[\big(1 + a(u)\big)M_+v\big]_r
  \Big(s, r_-\big[s; \beta\big(t_+(\sigma;
  \alpha), \sigma\big)\big]\Big),\\[-2mm]\\
\quad\quad \times\ \frac{d}{dr}r_-\Big[s;
  \beta\Big(t_+\big(\sigma; \alpha(t, r)\big),
  \sigma\Big)\Big]dsd\sigma,\\[-2mm]\\
A_{12} = - \frac{1}{r}\int^r_0\frac{1}{r}\int^t_{t_+(\sigma;
  \alpha)}\frac{a^\prime(u)(t - r) rM_-uu_t}{a(u)}\big(s,
  r_-(s; \beta)\big)dsd\sigma,\\[-2mm]\\
A_{13} = - \frac{1}{r}\int_0^r\frac{1}{r}\Big\{\int_0^{t_+(\sigma;
  \alpha)} \frac{a^\prime(u)(t - r)rM_-uu_t}{a(u)}\big(s, r_-(s;
  \beta)\big)ds,\\[-2mm]\\
\quad\quad - \int_0^{t_+(\sigma;
  \alpha)}\frac{a^\prime(u)(t - r) rM_-uu_t}{a(u)}\Big(s,
  r_-\big[s; \beta\big(t_+(\sigma; \alpha),
  \sigma\big)\big]\Big)ds\Big\}d\sigma,\\[-2mm]\\
A_{14} = - \frac{1}{r}\int_0^r\frac{a^\prime(u)(t - r)
  rM_-uu_t}{a(u)}\big(t_+(\sigma; \alpha), \sigma\big)\frac{d}{dr}
  t_+\big(\sigma; \alpha(t, r)\big)d\sigma,\\[-2mm]\\
A_{15} = - \frac{1}{r}\int_0^r\int_0^{t_+(\sigma;
  \alpha)}\Big[\frac{a^\prime(u)(t - r)rM_-uu_t}{a(u)}\Big]_r\Big(s,
  r_-\big[s; \beta\big(t_+(\sigma; \alpha),
  \sigma\big)\big]\Big),\\[-2mm]\\
\quad\quad \times\ \frac{d}{dr}r_-\Big[s;
  \beta\Big(t_+\big(\sigma; \alpha(t, r)\big),
  \sigma\Big)\Big]dsd\sigma.
\end{cases}
\end{equation}

Similarly, differentiating the equality \eqref{f5} with respect to
$t$, we have
\begin{eqnarray}\label{f8}
&&[(t - r)u_t](t, r) = B_1 + B_3 + B_6 + B_7 + B_{11} + B_{14} +
B_{15},
\end{eqnarray}
where
\begin{equation}\label{f9}
\begin{cases}
B_1 = - u(t, r),\\[-2mm]\\
B_3 = \frac{1}{r}\int_0^r\big(M_+[(t - r)v]\big)_r\big[0,
  \beta\big(t_+(\sigma; \alpha),\sigma\big)\big]
  \frac{d}{dt}\beta\Big(t_+\big(\sigma; \alpha(t, r)\big),
  \sigma\Big)d\sigma,\\[-2mm]\\
B_6 = \frac{1}{r}\int_0^r\Big[\frac{\big(1 - a(u)\big)v}
  {a(u)}\Big]_t\big(t_+(\sigma; \alpha), \sigma\big)\frac{
  dt_+\big(\sigma; \alpha(t, r)\big)}{dt}d\sigma,\\[-2mm]\\
B_7 =  - \frac{1}{r}\int_0^r\Big[\frac{\big(1 - a(u)\big)v}
  {a(u)}\Big]_r\Big(0, \beta\big(t_+(\sigma; \alpha),
  \sigma\big)\Big)\frac{d\beta\Big(t_+\big(\sigma; \alpha(t, r)\big),
  \sigma\Big)}{dt}d\sigma,\\[-2mm]\\
B_{11} = \frac{1}{r}\int_0^r\int_0^{t_+(\sigma;
  \alpha)}\big[\big(1 + a(u)\big)M_+v\big]_r
  \Big(s, r_-\big[s; \beta\big(t_+(\sigma;
  \alpha), \sigma\big)\big]\Big),\\[-2mm]\\
\quad\quad \times\ \frac{d}{dt}r_-\Big[s;
  \beta\Big(t_+\big(\sigma; \alpha(t, r)\big),
  \sigma\Big)\Big]dsd\sigma,\\[-2mm]\\
B_{14} = - \frac{1}{r}\int_0^r\frac{a^\prime(u)(t - r)
  rM_-uu_t}{a(u)}\big(t_+(\sigma; \alpha), \sigma\big)\frac{d}{dt}
  t_+\big(\sigma; \alpha(t, r)\big)d\sigma,\\[-2mm]\\
B_{15} = - \frac{1}{r}\int_0^r\int_0^{t_+(\sigma;
  \alpha)}\Big[\frac{a^\prime(u)(t - r)rM_-uu_t}{a(u)}\Big]_r\Big(s,
  r_-\big[s; \beta\big(t_+(\sigma; \alpha),
  \sigma\big)\big]\Big),\\[-2mm]\\
\quad\quad \times\ \frac{d}{dt}r_-\Big[s;
  \beta\Big(t_+\big(\sigma; \alpha(t, r)\big),
  \sigma\Big)\Big]dsd\sigma.
\end{cases}
\end{equation}

\begin{rem}\label{rem61}
We point out here that the term $r$ in the integrands of the
quantities \eqref{f7} and \eqref{f9} is essential so that we can
use energy to estimate the second order derivatives of $u$.
Another ingredient is the weight $t - r$ in the left hand side of
the above equalities \eqref{f6} and \eqref{f8}, which provides us
an extra decay $t^{-1}$ for $u_t$ and $u_r$ inside the cone
$\mathcal{D}$ when the integrands in \eqref{f7} and \eqref{f9} are
integrated outside the cone $\mathcal{D}$, where $\alpha$ (also $t
- r$) will be cancelled by the weaker decay $(1 + \alpha)^{- 1 +
\frac{2}{K}}$ of $\partial u$ and $\partial v$ explored in
\eqref{b20} and \eqref{b19}. When the integrands in \eqref{f7} and
\eqref{f9} are integrated inside the cone $\mathcal{D}$, the
weight $t - r$ does not help any more. However, the decay of
$\partial u$ is faster inside the cone $\mathcal{D}$, and
$E_2\big[u(t)\big]$ is also decaying there, as has been revealed
in \eqref{b21} and \eqref{e11}. The combination of these two
advantages results in the feasibility of the estimation for
maximal functions for large time in step 2 and 3.
\end{rem}

\bigskip
\textit{Step 2. Estimates for $A_j$ and $B_j$ for large time
except for $A_8$, $A_{11}$, $A_{13}$, $A_{15}$and $B_{11}$,
$B_{15}$.}
\bigskip

Our goal is to show the following estimate
\begin{equation}\label{f10}
\Big(\int_0^T\big[(1 + t)^{1 - \frac{4}{K}}\sup_{r \leq
\frac{t}{4} + 1}|\partial u(t, r)|\big]^2dt\Big)^{\frac{1}{2}}
\leq \frac{C_6K^2\epsilon}{2}.
\end{equation}
Introduce the notation
\begin{equation}\nonumber
\|p(t, r)\|_{L^\infty_r\big(\mathcal{D}(t)\big)L^2_w([T_1, T_2])}
= \Big(\int_{T_1}^{T_2}\big[(1 + t)^{- \frac{4}{K}}\sup_{r \leq
\frac{t}{4} + 1}|p(t, r)|\big]^2dt\Big)^{\frac{1}{2}}.
\end{equation}
Then, to show \eqref{f10}, by using \eqref{f6}, \eqref{f8} and
Lemma \ref{lem34}, we find that it suffices to show
\begin{equation}\label{f11}
\|\partial u(t, r)\|_{L^\infty_r\big(\mathcal{D}(t)\big)L^2_w([0,
3])} \leq \frac{C_6K^2\epsilon}{4}
\end{equation}
and
\begin{equation}\label{f12}
\begin{cases}
\|A_j\|_{L^\infty_r\big(\mathcal{D}(t)\big)L^2_w([3, T])} \leq
\frac{C_6K^2\epsilon}{2 \times 88},\\
\|B_j\|_{L^\infty_r\big(\mathcal{D}(t)\big)L^2_w([3, T])} \leq
\frac{C_6K^2\epsilon}{2 \times 88}
\end{cases}
\end{equation}
for all $A_j$ and $B_j$ defined in \eqref{f7} and \eqref{f9}.

\begin{rem}\label{rem62}
For $t \geq 3$ and $r \leq \frac{t}{4} + 1$, \eqref{c4} implies
that the plus characteristic passing through the point $(t, r)$
intersects the $\tau$ axis at $(\alpha, 0)$ with $\alpha =
\alpha(t, r)$. Moreover, $\alpha(t, r) > 1$, $\beta\big(\alpha(t,
r), 0\big) > 1$. Then Lemma \ref{lem34} implies that
\begin{eqnarray}\nonumber
&&\Big(\int_3^T\big[(1 + t)^{1 - \frac{4}{K}}\sup_{r \leq
  \frac{t}{4} + 1}|\partial u(t,
  r)|\big]^2dt\Big)^{\frac{1}{2}}\\\nonumber
&&= \sum_{j}\Big(\int_3^T\Big[(1 + t)^{- \frac{4}{K}}\sup_{r
  \leq \frac{t}{4} + 1}\frac{(1 + t)(|A_j| + |B_j|)}{t -
  r}\Big]^2dt\Big)^{\frac{1}{2}}\\\nonumber
&&\leq 2\sum_{j}\Big(\int_3^T\big[(1 + t)^{-
  \frac{4}{K}}\sup_{r \leq \frac{t}{4} + 1}(|A_j|
  + |B_j|)\big]^2dt\Big)^{\frac{1}{2}}\\\nonumber
&&\leq 2 \times 22 \times \frac{C_6K^2\epsilon}{2 \times
  88} \leq \frac{C_6K^2\epsilon}{4}.
\end{eqnarray}
Consequently, to show \eqref{f10}, it suffices to show \eqref{f11}
and \eqref{f12}.
\end{rem}

In this step, we estimate the quantities in \eqref{f12} except for
$A_8$, $A_{11}$, $A_{13}$, $A_{15}$ and $B_{11}$, $B_{15}$. The
estimation of $A_8$, $A_{11}$, $A_{13}$, $A_{15}$ and $B_{11}$,
$B_{15}$ involves energy decay properties and estimation for
Hardy-Littlewood Maximal functions and will be done in step 3. We
will estimate the quantity in \eqref{f11} in step 4.

Now we begin to estimate these quantities one by one. By the first
inequality in \eqref{b18}, it is easy to see that
\begin{eqnarray}\label{f13}
&&\big\|(|A_1| + |B_1|)\big\|_{L^\infty_r\big(
 \mathcal{D}(t)\big) L^2_w([3, T])}\\\nonumber
&&= \Big(\int_3^T\big[(1 + t)^{- \frac{4}{K}}\sup_{r \leq
  \frac{t}{4} + 1}(|A_1| + |B_1|)\big]^2dt
  \Big)^{\frac{1}{2}}\\\nonumber
&&\leq \Big(\int_3^T\big[(1 + t)^{- \frac{4}{K}}
  \frac{C\epsilon}{(1 + t)^{\frac{3}{5}}}\big]^2dt\Big)^{\frac{1}{2}}\\\nonumber
&&\leq C\epsilon \leq \frac{C_6K^2\epsilon}{88}.
\end{eqnarray}

Since the initial data is supported in the unit ball (see
\eqref{a4}), Remark \ref{rem62} implies that $A_2 = A_3 = A_5 =
A_7 = B_3 = B_7 = 0$.

Noting that \eqref{b21} and the second inequality in \eqref{b18},
we have
\begin{eqnarray}\label{f14}
&&\int_\alpha^t|u\partial u| \big(s, r_+(s;
  \alpha)\big)ds\\\nonumber
&&\leq CK\epsilon(1 + \alpha)^{- 1 + \frac{2}{K}}(1 + \alpha)^{
  - 1 + \frac{4}{K}}\int_\alpha^t(1 + s)^{1 - \frac{4}{K}}
  \sup_{\rho \leq \frac{s}{4} + 1}|\partial u(s, \rho)|ds\\\nonumber
&&\leq CK\epsilon(1 + \alpha)^{- 2 + \frac{6}{K}}(1 +
  t)^{\frac{1}{2}}\Big(\int_\alpha^t(1 + s)^{1 - \frac{4}
  {K}}\sup_{\rho \leq \frac{s}{4} + 1}|\partial u(s, \rho)|
  ds\Big)^{\frac{1}{2}}\\\nonumber
&&\leq CK^3\epsilon^2(1 + \alpha)^{- 2 + \frac{6}{K}}(1 +
  t)^{\frac{1}{2}}.
\end{eqnarray}
Similarly, by the second inequalities in \eqref{c5} and
\eqref{b18}, we have
\begin{eqnarray}\label{f15}
&&\frac{1}{r}\int_\alpha^t|u|^2\big(s, r_+(s;
  \alpha)\big)ds\\\nonumber
&&\leq \frac{CK^2\epsilon^2(t - \alpha)}{r(1 +
  \alpha)^{2 - \frac{4}{K}}} \leq CK^2\epsilon^2(1 +
  \alpha)^{ - 2 + \frac{4}{K}}.
\end{eqnarray}
Thus, noting the second inequality in \eqref{c6}, we obtain
\begin{eqnarray}\label{f16}
&&|A_4| = \Big|\frac{1}{r}\int_0^r\frac{1}{r}
  \int_\sigma^rM_+\Big[\frac{\big(1 - a(u)\big)v}{a(u)}\Big]
  \big(t_+(\rho; \alpha), \rho\big)d\rho d\sigma\Big|\\\nonumber
&&\leq \frac{C}{r}\int_\alpha^t\big(|uM_+v| + |vM_+u|\big)
  \big(s, r_+(s; \alpha)\big)ds \\\nonumber
&&\leq C\int_\alpha^t|u\partial u| \big(s, r_+(s; \alpha)\big)ds
  + \frac{C}{r}\int_\alpha^t|u|^2\big(s, r_+(s; \alpha)\big)ds\\\nonumber
&&\leq CK^3\epsilon^2(1 + t)^{- \frac{3}{2} + \frac{6}{K}},
\end{eqnarray}
which gives
\begin{eqnarray}\label{f17}
\|A_4\|_{L^\infty_r\big(\mathcal{D}(t)\big)L^2_w([3,
  T])} \leq \frac{C_6K^2\epsilon}{2 \times
88}.
\end{eqnarray}

To estimate $A_6$ and $B_6$, for $r \leq \frac{t}{4} + 1$, we use
the second line of \eqref{c1} and the first line of \eqref{c2} to
compute
\begin{eqnarray}\label{f18}
&&\Big|\frac{d}{dr}t_+\big(\sigma; \alpha(t,
  r)\big)d\sigma\Big| + \Big|\frac{d}{dt}t_+\big(\sigma; \alpha(t,
  r)\big)d\sigma\Big|\\\nonumber
&&\leq C\exp\Big\{C\Big(\int_\alpha^t(1 + s)^{- 2 +
  \frac{8}{K}}ds\Big)^{\frac{1}{2}}\\\nonumber
&&\quad\times \Big(\int_\alpha^t
  \Big[(1 + s)^{1 - \frac{4}{K}}\sup_{\rho \leq
  \frac{s}{4} + 1}|\partial u(s, \rho)|\Big]^2ds
  \Big)^{\frac{1}{2}}\Big\}\\\nonumber
&&\leq C\exp\Big\{CK^2\epsilon(1 + t)^{- \frac{1}{2} +
  \frac{4}{K}}\Big\} \leq C.
\end{eqnarray}
Consequently, the similar arguments as in \eqref{f14} give
\begin{eqnarray}\label{f19}
&&\big\|(|A_6| + |B_6|)\big\|_{L^\infty_r\big(
 \mathcal{D}(t)\big) L^2_w([3, T])}\\\nonumber
&&\leq  C\Big\|\frac{1}{r}\int_0^rr|uu_t|\big(t_+(\sigma;
  \alpha), \sigma\big)d\sigma\Big\|_{L^\infty_r\big(
 \mathcal{D}(t)\big) L^2_w([3, T])} \\\nonumber
&&\leq C\Big\|\int_\alpha^t|uu_t|
  \big(t_+(\sigma; \alpha), \sigma\big)d\sigma\Big\|_{L^\infty_r\big(
 \mathcal{D}(t)\big) L^2_w([3, T])}\\\nonumber
&&\leq \frac{C_6K^2\epsilon}{88}.
\end{eqnarray}

By \eqref{f18} and the first inequality of \eqref{b19}, $A_9$ and
$A_{10}$ can be easily estimated as
\begin{eqnarray}\label{f20}
&&\big\|(|A_9| + |A_{10}|)\big\|_{L^\infty_r\big(
 \mathcal{D}(t)\big) L^2_w([3, T])}\\\nonumber
&&\leq C\epsilon\big\|(1 + t)^{- 1 +
  \frac{1}{K}}\big\|_{L^\infty_r\big(
 \mathcal{D}(t)\big) L^2_w([3, T])}\\\nonumber
&&\leq \frac{C_6K^2\epsilon}{88}.
\end{eqnarray}

Note that for $r \leq \frac{t}{4} + 1$, it follows from \eqref{c4}
that
\begin{equation}\label{f21}
\frac{r_-(\alpha; \beta)}{r} \leq 3.
\end{equation}
Thus, by the second inequality in \eqref{c6} and Lemma
\ref{lem35}, we have
\begin{eqnarray}\label{f22}
&&|A_{12}| \leq \frac{C}{r}\int^r_0\frac{r_-(\alpha; \beta)}{r}
  \int^t_\alpha\big|(t - r)M_-uu_t\big(s,
  r_-(s; \beta)\big)\big|dsd\sigma\\\nonumber
&&\leq C\int^t_\alpha\big|(t - r)M_-uu_t\big(s,
  r_-(s; \beta)\big)\big|ds\\\nonumber
&&= C\int^t_\alpha\chi_{\{\frac{s}{4} + 1 < r_-(s; \beta)
   \leq \frac{2s}{3} + 1 \}}(s)\big|(t - r)M_-uu_t
  \big(s, r_-(s; \beta)\big)\big|ds\\\nonumber
&&\quad +\ C\int^t_\alpha\chi_{\{r_-(s; \beta) \leq \frac{s}{4} +
  1\}}(s)\big|(t - r)M_-uu_t\big(s, r_-(s; \beta)\big)\big|ds\\\nonumber
&&\leq CK^2\epsilon^2\int^t_\alpha\frac{s - r_-(s; \beta)}
  {\big[1 + \alpha\big(s, r_-(s; \beta)\big)\big]^{2 -
  \frac{4}{K}}(1 + s)^2}ds\\\nonumber
&&\quad +\ C\int^t_\alpha s(1 + s)^{- 2 +
  \frac{8}{K}}(1 + s)^{2 - \frac{8}{K}}\sup_{\rho \leq
  \frac{s}{4} + 1}|\partial u(s, \rho)|^2ds\\\nonumber
&&\leq CK^2\epsilon^2(1 + t)^{- \frac{3}{2}}\Big(
  \int^t_\alpha\frac{1}{\big[1 + \alpha\big(s,
  r_-(s; \beta)\big)\big]^{2 - \frac{8}{K}}}
  ds\Big)^{\frac{1}{2}}\\\nonumber
&&\quad +\ CK^3\epsilon^2(1 + t)^{- 1 + \frac{8}{K}}\\\nonumber
&&\leq CK^2\epsilon^2(1 + t)^{- \frac{3}{2} + CK^2\epsilon} +
  CK^3\epsilon^2(1 + t)^{- 1 + \frac{8}{K}}\\\nonumber
&&\leq CK^3\epsilon^2(1 + t)^{- 1 + \frac{8}{K}},
\end{eqnarray}
where $\chi_{I}(s)$ is the characteristic function defined on the
interval $I \subseteq \mathbb{R}$ by
\begin{eqnarray}\nonumber
\chi_{I}(s) = 1\quad {\rm for}\quad s \in I,\quad \chi_{I}(s) =
0\quad {\rm for}\quad s \in I^c = \mathbb{R}\setminus I.
\end{eqnarray}
\eqref{f22} gives that
\begin{eqnarray}\label{f23}
\|A_{12}\|_{L^\infty_r\big(
  \mathcal{D}(t)\big) L^2_w([3, T])} \leq \frac{C_6K^2\epsilon}{2
\times 88}.
\end{eqnarray}

Similarly, by \eqref{f18}, we have
\begin{eqnarray}\label{f24}
&&\big\|(|A_{14}| + |B_{14}|)\big\|_{L^\infty_r\big(
  \mathcal{D}(t)\big) L^2_w([3, T])}\\\nonumber
&&\leq C\Big\|\int_\alpha^t\big|(t - r)
  M_-uu_t\big(\tau, r_+(\tau; \alpha)\big)\big|d\tau\Big\|_
  {L^\infty_r\big(\mathcal{D}(t)\big) L^2_w([3, T])}\\\nonumber
&&\leq \frac{C_6K^2\epsilon}{88}.
\end{eqnarray}

\bigskip
\textit{Step 3. Estimate for  $A_8$, $A_{11}$, $A_{13}$,
$A_{15}$and $B_{11}$, $B_{15}$.}
\bigskip

First of all, we compute
\begin{eqnarray}\label{f25}
&&|A_8| \leq \frac{C}{r}\int_0^r\int_0^{t_+(\sigma; \alpha)}
  \frac{1}{r}\int_{r_-\big[s; \beta\big(t_+(\sigma;
  \alpha), \sigma\big)\big]}^{r_-(s; \beta)}\\\nonumber
&&\quad \big[|u_r\partial v| + |M_+u| + |rM_+u_r|\big]
  (s, \rho)d\rho dsd\sigma\\\nonumber
&&\leq C\int_{t_1}^t\frac{1}{r}\int_{r_-\big(s;
  \beta(\alpha, 0)\big)}^{r_-(s; \beta)}
  |u_r\partial v|(s, \rho)d\rho ds\\\nonumber
&&\quad +\ C\int_{t_1}^t\frac{1}{r}\int_{r_-\big(s;
  \beta(\alpha, 0)\big)}^{r_-(s; \beta)}\big(|M_+u| +
  |rM_+u_r|\big)(s, \rho)d\rho ds,
\end{eqnarray}
where $t_1$ is defined in \eqref{c7} (see Figure 2), and we used
the simplified notation
\begin{equation}\label{f25-1}
r_-\big(s; \beta(\alpha, 0)\big) =
\begin{cases} r_-\big(s; \beta(\alpha,
0)\big)\quad {\rm for}\quad  t_1 \leq s \leq \alpha,\\
0\quad {\rm for}\quad  \alpha < s \leq t.
\end{cases}
\end{equation}
Similarly as in \eqref{f22} and using Lemma \ref{lem36}, we have
\begin{eqnarray}\label{f26}
&&\int_{t_1}^t\frac{1}{r}\int_{r_-\big(s;
  \beta(\alpha, 0)\big)}^{r_-(s; \beta)}
  |u_r\partial v(s, \rho)|d\rho ds\\\nonumber
&&\leq \int_{t_1}^t\frac{1}{r}\int_{r_-\big(s;
  \beta(\alpha, 0)\big)}^{r_-(s; \beta)}\frac{CK\epsilon^2\chi_{\{\rho
  \geq \frac{s}{4} + 1 \}}(s)}{(1 + s)\big[1 + \alpha(s, r_-(s;
  \beta))\big]^{1 - \frac{1}{K}}}d\rho ds\\\nonumber
&&\quad +\ \int_{t_1}^t\frac{1}{r}
  \int_{r_-\big(s; \beta(\alpha, 0)\big)}^{r_-(s; \beta)}
  \frac{C\epsilon(1 + s)^{1 - \frac{4}{K}}\sup_{\rho
  \leq \frac{s}{4} + 1 }|\partial u(s, \rho)|}
  {(1 + s)^{1 - \frac{1}{K}}(1
  + s)^{1 - \frac{4}{K}}}d\rho ds\\\nonumber
&&\leq CK\epsilon^2(1 + t)^{- \frac{1}{2} + CK^2\epsilon}
  \Big(\int_{t_1}^t\frac{1}{\big[1 + \alpha\big(s, r_-(s; \beta)
  \big)\big]^{2 - \frac{2}{K}}}ds\Big)^{\frac{1}{2}}\\\nonumber
&&\quad +\ CK^2\epsilon^2(1 + t)^{- \frac{3}{2} + \frac{5}{K} +
  CK^2\epsilon}\\\nonumber
&&\leq CK^2\epsilon^2(1 + t)^{- \frac{1}{2} + CK^2\epsilon},
\end{eqnarray}
and

\begin{figure}
\medskip
\includegraphics [width=12cm,clip]{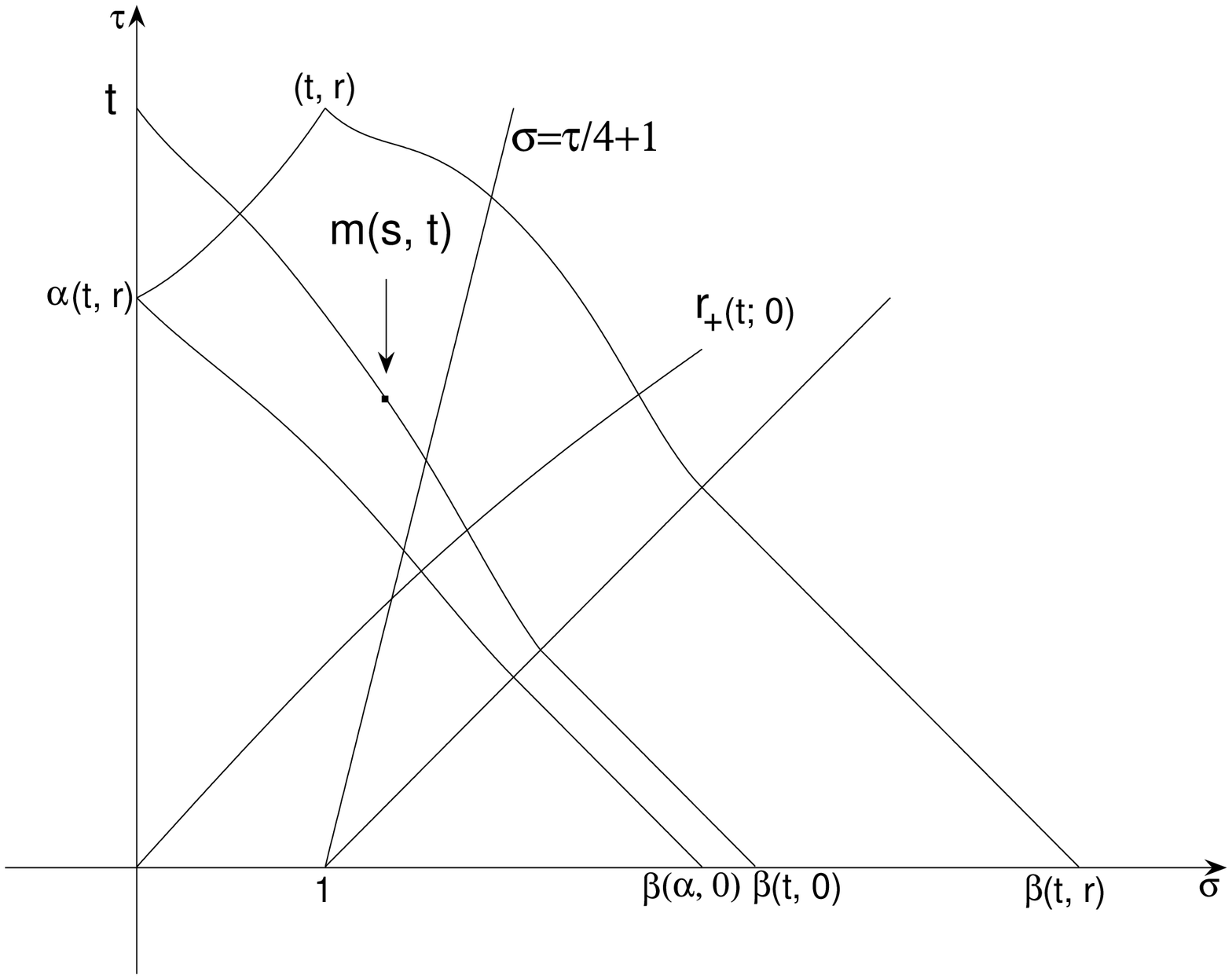}\\
\caption{}
\medskip
\end{figure}

\begin{eqnarray}\label{f27}
&&\int_{t_1}^t\frac{1}{r}\int_{r_-\big(s;
  \beta(\alpha, 0)\big)}^{r_-(s; \beta)}\big(|M_+u| +
  |rM_+u_r|\big)(s, \rho)d\rho ds\\\nonumber
&&\leq C(1 +
t)^{CK^2\epsilon}\int_{t_1}^t\mathbb{M}\big[\big(|M_+u|
  + |rM_+u_r|\big)(s, \rho)\big]\big(m(s, t)\big)ds\\\nonumber
&&= \int_{t_1}^tC(1 + s)^{CK^2\epsilon}\mathbb{M}\big[\big(|M_+u|
+
  |rM_+u_r|\big)(s, \rho)\big]\big(m(s, t)\big)ds,
\end{eqnarray}
where $m(s, t)$ is defined in Lemma \ref{lem37} (see Figure 3). We
point out here that in the above inequality \eqref{f26},
$\alpha\big(s, r_-(s; \beta)\big)$ may also be $\gamma\big(s,
r_-(s; \beta)\big) \in [0, 1]$. Since we can treat them exactly
the same way, we did not distinguish them, and in the following
context, we sometimes also do the same thing. Consequently, the
combination of \eqref{e10}, Lemma \ref{lem37} and Hardy inequality
imply that
\begin{eqnarray}\label{f28}
&&\|A_8\|_{L^\infty_r\big(\mathcal{D}(t)\big) L^2_w([3,
  T])}\\\nonumber
&&\leq CK^3\epsilon^2 + \Big(\int_3^T\Big[\int_{t_1}^t
  (1 + s)^{- \frac{4}{K} + CK^2
  \epsilon}\mathbb{M}\big[\big(|M_+u|\\\nonumber
&&\quad +\ |rM_+u_r|\big)(s, \rho)\big]\big(m(s,
  t)\big)ds\Big]^2dt\Big)^{\frac{1}{2}}\\\nonumber
&&\leq CK^3\epsilon^2 + \int_3^T(1 + s)^{- \frac{4}{K} +
  CK^2\epsilon}\big\|\mathbb{M}\big[\big(|M_+u|\\\nonumber
&&\quad +\ |rM_+u_r|\big)(s, \rho)\big)\big(m(s,
  t)\big)\big\|_{L^2(0 \leq t \leq T)}ds\\\nonumber
&&\leq CK^3\epsilon^2 + \int_3^T(1 + s)^{- \frac{4}{K} +
  CK^2\epsilon}\Big(\Big\|r\frac{M_+u}{|x|}(s,
  \rho)\Big\|_{L^2}\\\nonumber
&&\quad +\ \big\|rM_+u_r(s, \rho)\big\|_{L^2}\Big)ds\\\nonumber
&&\leq CK^3\epsilon^2 + \int_3^T(1 + s)^{- \frac{4}{K} +
  CK^2\epsilon}\Big(\Big\|M_+u_r(s, x)\Big\|_{L^2(R^3)}\\\nonumber
&&\quad +\ \|u_r^2(s, x)\|_{L^2(R^3)}\Big)ds\\\nonumber&&\leq
  CK^2\epsilon^2 + \int_3^T(1 + s)^{- \frac{4}{K} +
  CK^2\epsilon}\Big(\frac{\|\widetilde{\Gamma}_1u_r(s, x)
  \|_{L^2(R^3)}}{1 + s}\\\nonumber
&&\quad +\ \big[\|u(s, \cdot)\|_{L^\infty} + \|\partial u(s,
  \cdot)\|_{L^\infty}\big]\|u_r(s, x)\|_{L^2(R^3)}\Big)ds\\\nonumber
&&\leq CK^3\epsilon^2 + \int_3^TCE_2\big(u(0)\big)^{\frac{1}{2}}(1
  + s)^{- \frac{4}{K} + CK^2\epsilon}\Big[(1 + s)^{- 1 +
  \frac{3}{K}}\\\nonumber
&&\quad +\ CK\epsilon(1 + s)^{- 1 + \frac{1}{K}} + \|\partial
  u(s, \cdot)\|_{L^\infty}(1 + s)^{CK^2\epsilon}\Big]ds\\\nonumber
&&\leq CKE_2\big(u(0)\big)^{\frac{1}{2}} + \frac{C_6K^2\epsilon}{4
\times 88}.
\end{eqnarray}

Similarly as in \eqref{f18}, we use the fifth equality in
\eqref{c1} and the third line in \eqref{c2} to calculate
\begin{eqnarray}\label{f29}
&&\Big|\frac{d}{dr}r_-\Big[s;
  \beta\Big(t_+\big(\sigma; \alpha(t, r)\big),
  \sigma\Big)\Big]\Big|\\\nonumber
&&\quad +\ \Big|\frac{d}{dt}r_-\Big[s;
  \beta\Big(t_+\big(\sigma; \alpha(t, r)\big),
  \sigma\Big)\Big]\Big|\\\nonumber
&&\leq C\Big|r_{-, \beta}\Big[s;
  \beta\Big(t_+\big(\sigma; \alpha(t, r)\big),
  \sigma\Big)\Big]\beta_t\Big(t_+\big(\sigma; \alpha(t, r)\big),
  \sigma\Big)\Big|\\\nonumber
&&\leq C\exp\Big\{C\int_s^{t_+(\sigma; \alpha)}
  \Big|u_r\Big[\tau, r_-\Big(\tau; \beta\big[t_+\big(\sigma; \alpha(t, r)\big),
  \sigma\big]\Big)\Big]\Big|d\tau\\\nonumber
&&\leq C\exp\Big\{\int_s^{s_3}\frac{CK\epsilon}{1 + \tau}d\tau +
  C\int_{s_3}^t\sup_{\sigma \leq \frac{\tau}{4} + 1}
  |\partial u(\tau, \sigma)|d\tau\Big\}\\\nonumber
&&\leq C\exp\Big\{\int_s^{s_3}\frac{CK\epsilon}{1 + \tau}d\tau +
  C\int_{s_3}^t\frac{(1 + s)^{1 - \frac{4}{K}}\sup_{\sigma \leq \frac{\tau}{4} + 1}
  |\partial u(\tau, \sigma)|}{(1 + s)^{1 - \frac{4}{K}}}d\tau\Big\}\\\nonumber
&&\leq C,
\end{eqnarray}
where $s_3$ is determined by (see Figure 2 as reference)
\begin{eqnarray}\nonumber
r_-\Big(s_3; \beta\big[t_+\big(\sigma; \alpha(t, r)\big),
  \sigma\big]\Big) = \frac{s_3}{4} + 1.
\end{eqnarray}

Thus, similarly as in \eqref{f25}-\eqref{f28}
\begin{eqnarray}\nonumber
|A_{11}| + |B_{11}| \leq
  \frac{1}{r}\int_0^r\int_0^t\big[|\partial u\partial v| + r|M_+u_r|\big]
  \Big(s, r_-\big[s; \beta\big(t_+(\sigma;
  \alpha), \sigma\big)\big]\Big)dsd\sigma,
\end{eqnarray}
which gives
\begin{eqnarray}\label{f30}
\big\|(|A_{11}| + |B_{11}|)\big\|_{L^\infty_r\big(
  \mathcal{D}(t)\big) L^2_w([3, T])} \leq
CKE_2\big(u(0)\big)^{\frac{1}{2}} + \frac{C_6K^2\epsilon}{2 \times
88}.
\end{eqnarray}

Similarly as in \eqref{f25}, we calculate
\begin{eqnarray}\label{f31}
&&|A_{13}| \leq \frac{1}{r}\int_0^r\int_0^t\frac{1}{r}\int^{r_-(s;
  \beta)}_{r_-\big(s; \beta(\alpha, 0)\big)}\Big[
  \frac{a^\prime(u)(t - r)rM_-uu_t}{a(u)}\Big]_r(s,
  \rho)d\rho dsd\sigma\\\nonumber
&&\leq C\int_{t_1}^t\frac{1}{r}\int^{r_-(s;
  \beta)}_{r_-\big(s; \beta(\alpha, 0)\big)}\big(
  |u_r(t - r)rM_-uu_t| + |tM_-uu_t|\big)(s,
  \rho)d\rho ds\\\nonumber
&&\quad +\ C(1 + t)^{CK^2\epsilon}\int_{t_1}^t\mathbb{M}\Big(\big[
  (t - r)r(M_-uu_t)_r\big](s, \rho)\Big)\big(m(s, t)\big) ds,
\end{eqnarray}
where the meaning of $r_-\big(s; \beta(\alpha, 0)\big)$ is the
same as \eqref{f25-1}. By \eqref{b19}, \eqref{b20}, \eqref{b21}
and Lemma \ref{lem35}, we have
\begin{eqnarray}\label{f32}
&&\int_{t_1}^t\frac{1}{r}\int^{r_-(s;
  \beta)}_{r_-\big(s; \beta(\alpha, 0)\big)}\big(
  |u_r(t - r)rM_-uu_t| + |tM_-uu_t|\big)(s,
  \rho)d\rho ds\\\nonumber
&&\leq \int_{t_1}^t
  \frac{1}{r} \int^{r_-(s; \beta)}_{r_-\big(s; \beta(\alpha,
  0)\big)}\Big\{\frac{CK^3\epsilon^3s(s-\rho)}{\big(1 + \alpha(s, \beta)\big)^{3
  - \frac{6}{K}}}\frac{\chi_{\{\rho \geq \frac{s}{4}
  + 1\}}(s)}{(1 + s)^3}\\\nonumber
&&\quad +\ \frac{C\epsilon s(1 + s)^{- 2 + \frac{8}{K}}}{(1 +
  s)^{1 - \frac{1}{K}}}\Big[(1 + s)^{- 2 +
  \frac{8}{K}}|\partial u(s, \rho)|^2\Big]\chi_{\{\rho \leq
  \frac{s}{4} + 1\}}(s)\\\nonumber
&&\quad +\ C s(1 + s)^{- 2 + \frac{8}{K}}\Big[(1 + s)^{2 -
  \frac{8}{K}}|\partial u(s, \rho)|^2\Big]\chi_{\{\rho \leq
  \frac{s}{4} + 1\}}(s)\\\nonumber
&&\quad +\ \frac{CK^2\epsilon^2s}{(1 + s)^2}\frac{\chi_{\{\rho
  \geq \frac{s}{4} + 1\}}(s)}{(1 + \alpha)^{2
  - \frac{4}{K}}}\Big\}d\rho ds\\\nonumber
&&\leq CK^4\epsilon^2(1 + t)^{- 1 + CK^2\epsilon}.
\end{eqnarray}
\eqref{e11} and similar arguments as in \eqref{f25}-\eqref{f28}
give
\begin{eqnarray}\label{f33}
&&\Big(\int_3^T\Big[\int_{t_1}^t(1 + s)^{- \frac{4}{K} +
  CK^2\epsilon}\mathbb{M}\Big(\big[(t - r)r\partial u\partial^2u
  \big](s, \rho)\Big)\big(m(s, t)\big)ds\Big]^2dt
  \Big)^{\frac{1}{2}}\\\nonumber
&&\leq \int_3^T(1 + s)^{- \frac{4}{K} + CK^2\epsilon}
  \big\|\big[(t - r)\partial u\partial^2u\big](s,
  x)\big\|_{L^2(\mathbb{R}^3)}ds\\\nonumber
&&\leq \int_3^T(1 + s)^{- \frac{4}{K} + CK^2\epsilon}
  \Big[\\\nonumber
&&\quad  \big\|(s - r)\partial u(s, x)\big\|_{L^\infty(|x|
  \geq \frac{s}{4} + 1)}\big\|\partial^2u(s, x)(s, x)
  \big\|_{L^2(\mathbb{R}^3)}\\\nonumber
&&\quad +\ \big\|\partial u(s, x)\big\|_{L^\infty(|x|
  \leq \frac{s}{4} + 1)}\big\|\widetilde{\Gamma}\partial u(s, x)(s, x)
  \big\|_{L^2(|x| \leq \frac{s}{4} + 1)}\Big]ds\\\nonumber
&&\leq \int_3^T(1 + s)^{- \frac{4}{K} + CK^2\epsilon}
  \Big[\frac{CK\epsilon E_2\big(u(0)\big)^{\frac{1}{2}}(1 +
  s)^{CK^2\epsilon}(s - r)}{(1 + s)\big(1 + \alpha(s, r)\big)^{1
  - \frac{2}{K}}}\\\nonumber
&&\quad +\ CE_2\big(u(0)\big)^{\frac{1}{2}}(1 +
  s)^{- 1 + \frac{3}{K}}\big\|\partial u(s, x)\big\|_{L^\infty(|x|
  \leq \frac{s}{4} + 1)}\Big]ds\\\nonumber
&&\leq CK\epsilon E_2\big(u(0)\big)^{\frac{1}{2}}.
\end{eqnarray}
Combining \eqref{f31}-\eqref{f33}, we have
\begin{eqnarray}\label{f34}
\|A_{13}\|_{L^\infty_r\big(
  \mathcal{D}(t)\big) L^2_w([3, T])}  \leq
CK\epsilon E_2\big(u(0)\big)^{\frac{1}{2}} +
\frac{C_6K^2\epsilon}{4 \times 88}.
\end{eqnarray}

At last, we  estimate $A_{15}$ and $B_{15}$ as
\begin{eqnarray}\label{f35}
&&|A_{15}| + |B_{15}|\\\nonumber &&\leq \frac{C}{r}\int_0^r
  \int_0^t\Big|(t - r)r(M_-uu_t)_r\Big(s, r_-\big[s;
  \beta\big(t_+(\sigma; \alpha), \sigma\big)\big]
  \Big)\Big|dsd\sigma\\\nonumber
&&\quad +\ \frac{C}{r}\int_0^r\int_0^t
  \big[|u_r(t - r)rM_-uu_t|\\\nonumber
&&\quad +\ |tM_-uu_t|\big]\Big(s,
  r_-\big[s; \beta\big(t_+(\sigma; \alpha),
  \sigma\big)\big]\Big)dsd\sigma.
\end{eqnarray}
Similarly, we have
\begin{eqnarray}\label{f36}
\big\|(|A_{15}| + |B_{15}|)\big\|_{L^\infty_r\big(
  \mathcal{D}(t)\big) L^2_w([3, T])} \leq CK\epsilon E_2\big(u(0)\big)^{\frac{1}{2}} +
\frac{C_6K^2\epsilon}{2 \times 88}.
\end{eqnarray}

\bigskip
\textit{Step 4. Estimate for the quantity in \eqref{f11}.}
\bigskip

In this case, it is obviously that it suffices to show
\begin{equation}\label{f37}
\Big(\int_0^3\sup_{r \leq \frac{t}{4} + 1}|\partial u(t,
r)|\big]^2dt\Big)^{\frac{1}{2}} \leq \frac{C_6K^2\epsilon}{8}.
\end{equation}

Similarly as in \eqref{f5}, for $t \geq t_+(r; 0)$, we integrate
equation \eqref{f1} to get
\begin{eqnarray}\label{f37-1}
&&u(t, r) = \frac{1}{r}\int_0^rM_+v\big[0, \beta\big(t_+(\sigma;
  \alpha), \sigma\big)\big]d\sigma\\\nonumber
&&\quad -\ \frac{1}{r}\int_0^r\int_0^{t_+(\sigma; \alpha)}
  \frac{a^\prime(u)rM_-uu_t}{a(u)}\Big(s, r_-\big[s;
  \beta\big(t_+(\sigma; \alpha), \sigma\big)\big]\Big)dsd\sigma,
\end{eqnarray}
Similarly as in \eqref{f6} and \eqref{f8}, differentiating the
above equation with respect to $t$ and $r$ yields
\begin{eqnarray}\\\nonumber
&&u_r(t, r) = \frac{1}{r}\int_0^r\frac{M_+v(0, \beta) -
  M_+v\big[0, \beta\big(t_+(\sigma; \alpha),\sigma\big)
  \big]}{r}d\sigma\\\nonumber
&&\quad +\ \frac{1}{r}\int_0^r(M_+v)_r\big[0,
  \beta\big(t_+(\sigma; \alpha),\sigma\big)\big]
  \frac{d}{dr}\beta\Big(t_+\big(\sigma; \alpha(t,
  r)\big),\sigma\Big)d\sigma\\\nonumber
&&\quad -\ \frac{1}{r}\int_0^r\int_0^{t_+(\sigma;
  \alpha)}\Big\{\frac{a^\prime(u)
  rM_-uu_t}{a(u)}\big(s, r_-(s; \beta)\big)\\\nonumber
&&\quad -\ \frac{a^\prime(u)rM_-uu_t}{a(u)}\Big(s,
  r_-\big[s; \beta\big(t_+(\sigma; \alpha),
  \sigma\big)\big]\Big)\Big\}dsd\sigma\\\nonumber
&&\quad -\ \frac{1}{r}\int_0^r\frac{1}{r}\int_{t_+(\sigma;
  \alpha)}^t \frac{a^\prime(u)
  rM_-uu_t}{a(u)}\big(s, r_-(s; \beta)\big)dsd\sigma\\\nonumber
&&\quad -\ \frac{1}{r}\int_0^r\frac{a^\prime(u)rM_-uu_t}
  {a(u)}\big(t_+(\sigma; \alpha), \sigma\big)\frac{d}{dr}
  t_+\big(\sigma; \alpha(t, r)\big)d\sigma\\\nonumber
&&\quad -\ \frac{1}{r}\int_0^r\int_0^{t_+(\sigma;
  \alpha)}\Big[\frac{a^\prime(u)rM_-uu_t}{a(u)}\Big]_r\Big(s,
  r_-\big[s; \beta\big(t_+(\sigma; \alpha),
  \sigma\big)\big]\Big)\\\nonumber
&&\quad \times\ \frac{d}{dr}r_-\Big[s;
  \beta\Big(t_+\big(\sigma; \alpha(t, r)\big), \sigma\Big)\Big]dsd\sigma.
\end{eqnarray}
and
\begin{eqnarray}\label{f37-3}
&&u_t(t, r) = \frac{1}{r}\int_0^r(M_+v)_r\big[0,
  \beta\big(t_+(\sigma; \alpha),\sigma\big)\big]
  \frac{d}{dt}\beta\Big(t_+\big(\sigma; \alpha(t,
  r)\big),\sigma\Big)d\sigma\\\nonumber
&&\quad -\ \frac{1}{r}\int_0^r\frac{a^\prime(u)rM_-uu_t}
  {a(u)}\big(t_+(\sigma; \alpha), \sigma\big)\frac{d}{dt}
  t_+\big(\sigma; \alpha(t, r)\big)d\sigma\\\nonumber
&&\quad -\ \frac{1}{r}\int_0^r\int_0^{t_+(\sigma;
  \alpha)}\Big[\frac{a^\prime(u)rM_-uu_t}{a(u)}\Big]_r\Big(s,
  r_-\big[s; \beta\big(t_+(\sigma; \alpha),
  \sigma\big)\big]\Big)\\\nonumber
&&\quad \times\ \frac{d}{dt}r_-\Big[s; \beta\Big(t_+\big(\sigma;
  \alpha(t, r)\big), \sigma\Big)\Big]dsd\sigma.
\end{eqnarray}
Thus, we have
\begin{eqnarray}\nonumber
&&|\partial u|(t, r) \leq \frac{1}{r}\int_0^r\Big|\frac{M_+v(0,
  \beta) - M_+v\big[0, \beta\big(t_+(\sigma; \alpha),
  \sigma\big)\big]}{r}\Big|d\sigma\\\nonumber
&&\quad +\ \frac{C}{r}\int_0^r\big|(M_+v)_r\big[0,
  \beta\big(t_+(\sigma; \alpha),\sigma\big)\big]\big|d\sigma\\\nonumber
&&\quad +\ \frac{C}{r}\int_0^r\int_0^{t_+(\sigma;
  \alpha)}\frac{1}{r}\int_{r_-\big(s; \beta(\alpha, 0)\big)}^{r_-(s; \beta)}
  \Big|\Big\{\frac{a^\prime(u)rM_-uu_t}{a(u)}
  \Big\}_r(s, \rho)\Big|d\rho dsd\sigma\\\nonumber
&&\quad +\ \frac{1}{r}\int_\alpha^t\big|rM_-uu_t\big(s, r_-(s;
  \beta)\big)\big|ds + \frac{C}{r}\int_0^r\big|rM_-uu_t\big(t_+(\sigma;
  \alpha), \sigma\big)\big|d\sigma\\\nonumber
&&\quad +\ \frac{C}{r}\int_0^r\int_0^{t_+(\sigma;
  \alpha)}\Big|\Big[\frac{a^\prime(u)rM_-uu_t}{a(u)}\Big]_r\Big(s,
  r_-\big[s; \beta\big(t_+(\sigma; \alpha),
  \sigma\big)\big]\Big)\Big|dsd\sigma.
\end{eqnarray}

Since $t \leq 3$, there holds
\begin{eqnarray}\nonumber
&&\frac{1}{r}\int_0^r\Big|\frac{M_+v(0,
  \beta) - M_+v\big[0, \beta\big(t_+(\sigma; \alpha),
  \sigma\big)\big]}{r}\Big|d\sigma\\\nonumber
&&\quad +\ \frac{C}{r}\int_0^r\big|(M_+v)_r\big[0,
  \beta\big(t_+(\sigma; \alpha),\sigma\big)\big]\big|d\sigma\\\nonumber
&&\leq \frac{1}{r}\int_0^r\frac{1}{r}\int_{\beta\big(t_+(\sigma;
  \alpha), \sigma\big)}^{\beta}|(M_+v)_r(0, \rho)d\rho|d\sigma\\\nonumber
&&\quad +\ \frac{C}{r}\int_{\beta(\alpha, 0)}^\beta|(M_+v)_r(0,
  \rho)|\frac{1}{\Big|\frac{d\beta\big(t_+(\sigma;
  \alpha), \sigma\big)}{d\sigma}\Big|}d\rho\\\nonumber
&&\leq \frac{C}{r}\int_{\beta(\alpha, 0)}^\beta|(M_+v)_r(0,
  \rho)|d\rho\\\nonumber
&&\leq C\mathbb{M}\big((M_+v)_r(0, \rho)\big)\big(\beta(t,
  0)\big).
\end{eqnarray}
Similar arguments give
\begin{eqnarray}\nonumber
&&\frac{1}{r}\int_0^r\int_0^{t_+(\sigma;
  \alpha)}\frac{1}{r}\int_{r_-\big(s; \beta(\alpha, 0)\big)}^{r_-(s; \beta)}
  \Big|\Big\{\frac{a^\prime(u)rM_-uu_t}{a(u)}
  \Big\}_r(s, \rho)\Big|d\rho dsd\sigma\\\nonumber
&&\leq \int_0^3\mathbb{M}\Big\{\Big[\frac{a^\prime(u)
  rM_-uu_t}{a(u)}\Big]_r(s, \rho)\Big\}\big(m(s, t)\big)ds.
\end{eqnarray}
Recall \eqref{f21}, we have
\begin{eqnarray}\nonumber
&&\frac{1}{r}\int_\alpha^t
  \big|rM_-uu_t\big(s, r_-(s; \beta)\big)\big|ds\\\nonumber
&&\quad +\ \frac{C}{r}\int_0^r\big|rM_-uu_t\big(t_+(\sigma;
  \alpha), \sigma\big)\big|d\sigma\\\nonumber
&&\leq C\int_\alpha^t\Big(\big|M_-uu_t\big(s, r_-(s;
  \beta)\big)\big| + \big|rM_-uu_t\big(t_+(\sigma;
  \alpha), \sigma\big)\big|\Big)ds\\\nonumber
&&\leq \int_0^3\|\partial u(s, \cdot)\|_{L^\infty}^2ds \leq
  CK^4\epsilon^2.
\end{eqnarray}
Variable of integration gives
\begin{eqnarray}\nonumber
&&\frac{1}{r}\int_0^r\int_0^{t_+(\sigma;
  \alpha)}\Big|\Big[\frac{a^\prime(u)rM_-uu_t}{a(u)}\Big]_r\Big(s,
  r_-\big[s; \beta\big(t_+(\sigma; \alpha),
  \sigma\big)\big]\Big)\Big|dsd\sigma\\\nonumber
&&\leq \int_0^t \frac{1}{r}\int_0^r\Big|\Big[\frac{a^\prime(u)
  rM_-uu_t}{a(u)}\Big]_r\Big(s, r_-\big[s; \beta\big(t_+(\sigma; \alpha),
  \sigma\big)\big]\Big)\Big|d\sigma ds\\\nonumber
&&\leq \int_0^t \frac{1}{r}\int_{r_-\big(s; \beta(\alpha, 0)
  \big)}^{r_-(s; \beta)}\Big|\Big[\frac{a^\prime(u)
  rM_-uu_t}{a(u)}\Big]_r(s, \rho)\Big|d\rho ds\\\nonumber
&&\leq \int_0^3\mathbb{M}\Big\{\Big[\frac{a^\prime(u)
  rM_-uu_t}{a(u)}\Big]_r(s, \rho)\Big\}\big(m(s, t)\big)ds.
\end{eqnarray}
Finally, we arrive at
\begin{eqnarray}\label{f38}
&&\big\|\partial u(t, r)\big\|_{L^\infty_r\big(
  \mathcal{D}(t)\big) L^2_w([0, 3])}
   \leq CK^4\epsilon^2\\\nonumber
&&\quad +\ C\big\|\mathbb{M}\big((M_+v)_r(0, \rho)\big)
  \big(\beta(t, 0)\big)\big\|_{L_t^2([0, 3])}\\\nonumber
&&\quad +\ \int_0^3\Big\|\mathbb{M}\Big\{\Big[\frac{a^\prime(u)
  rM_-uu_t}{a(u)}\Big]_r(s, \rho)\Big\}
  \big(m(s, t)\big)\Big\|_{L^2_t([0, 3])}ds\\\nonumber
&&\leq CK^4\epsilon^2 + C\big\|(M_+v)_r(0,
  r)\big\|_{L^2}\\\nonumber &&\quad +\
\int_0^3\Big\|\Big[\frac{a^\prime(u)
  rM_-uu_t}{a(u)}\Big]_r(s, r)\Big\|_{L^2}ds\\\nonumber
&&\leq CK^4\epsilon^2 + C\big\|\partial u(0, r)\big\|_{L^2} +
  C\big\|r\partial^2 u(0, r)\big\|_{L^2} + C\big\|r(\partial u)^2(0,
  r)\big\|_{L^2}\\\nonumber
&&\quad +\ C\int_0^3\Big(\big\|r(\partial u)^3(s, r)\big\|_{L^2} +
  \big\|r\partial u\partial^2 u(0, r)\big\|_{L^2} + \big\|(\partial u)^2(s,
  r)\big\|_{L^2}\Big)ds\\\nonumber
&&\leq CK^4\epsilon^2 + C\Big\|\frac{\partial u(0, x)}{|x|}
  \Big\|_{L^2(\mathbb{R}^3)} + CE_2\big[u(0)\big]^{\frac{1}{2}}
  + C\big\|\partial u(0, \cdot)\big\|_{L^4(\mathbb{R}^3)}^2\\\nonumber
&&\quad +\ C\int_0^3\Big(\big\|\partial u(s,
  \cdot)\big\|_{L^6(\mathbb{R}^3)}^3 +
  \big\|\partial u(s, \cdot)\big\|_{L^\infty}E_2\big[u(s)
  \big]^{\frac{1}{2}}\\\nonumber
&&\quad +\ CK^4\epsilon^2 + \sup_{\rho \leq \frac{s}{4} +
  1}|\partial u(s, \rho)|^2\Big)
  ds\\\nonumber
&&\leq CK^4\epsilon^2 + CE_2\big[u(0)\big]^{\frac{1}{2}}.
\end{eqnarray}

For $t \geq t_+(r; 0)$, similarly as in \eqref{f37-1}, we have
\begin{eqnarray}\label{f42}
&&u(t, r) = \frac{1}{r}\int_\gamma^rM_+v\big[0, \beta\big(
  t_+(\sigma; \gamma), \sigma\big)\big]d\sigma\\\nonumber
&&\quad -\ \frac{1}{r}\int_\gamma^r\int_0^{t_+(\sigma; \gamma)}
  \frac{a^\prime(u)rM_-uu_t}{a(u)}\Big(s, r_-\big[s;
  \beta\big(t_+(\sigma; \gamma), \sigma\big)\big]\Big)dsd\sigma,
\end{eqnarray}
By the same procedure above, we can prove
\begin{eqnarray}\label{f43}
\big\|\partial u(t, r)\big\|_{L^\infty_r\big(
  \mathcal{D}(t)\big) L^2_w([0, 3])}
   \leq CK^4\epsilon^2 + CE_2\big[u(0)\big]^{\frac{1}{2}}.
\end{eqnarray}

Finally, combining \eqref{f13}, \eqref{f17},
\eqref{f19},\eqref{f20}, \eqref{f23}, \eqref{f24}. \eqref{f28},
\eqref{f30}, \eqref{f34}, \eqref{f36} and \eqref{f38}, \eqref{f43}
all together, we have proved \eqref{b21} with $C_6$ being replaced
by $\frac{C_6}{2}$ provided that \eqref{a5} is satisfied.

\section{Stability of Local Strong Solutions}

Former sections showed that the Cauchy problem of the quasi-linear
wave equation \eqref{a1}-\eqref{a2} admits a global strong
solution under the hypothesis $(H_1)-(H_2)$. In this section, we
first show the local existence of general (may not be small)
strong solutions to \eqref{a1}-\eqref{a2} under the assumption
\begin{equation}\label{g0}
(\overline{H_1}):
\begin{cases} {\rm f\ and\ g\ are\ radially\
symmetric\ functions},\\
{\rm supp\big\{(f, g)\big\}\ is\ compact,\ }
\end{cases}
\end{equation}
then prove the stability of the local strong solutions by using
the homotopy method. The uniqueness of the global strong solution
under the hypothesis $(H_1)-(H_2)$ follows as a direct corollary
of the stability for local strong solutions under the hypothesis
$(\overline{H_1})$, which completes the proof of Theorem
\ref{thm11}. Moreover, as a byproduct, we in fact prove that the
local well-posedness of the Cauchy problem \eqref{a1}-\eqref{a2}
holds in the case of radially symmetric initial data $(f, g) \in
H^2(\mathbb{R}^3) \times H^1(\mathbb{R}^3)$ (the sharp local
well-posedness results for the general quasi-linear wave equations
holds when the initial data are in $H^s(\mathbb{R}^3) \times H^{s
- 1}(\mathbb{R}^3)$ for $s > 2$, see \cite{SmithTataru}).

As what we have discussed in section 2, to prove the local
existence of the quasi-linear wave equation \eqref{a1}-\eqref{a2},
we only need to show \textit{a priori} bound (depending only on
$\big\|(f, g)\big\|_{H^2(\mathbb{R}^3) \times H^1(\mathbb{R}^3)}$)
of the quantity
\begin{equation}\label{g1}
\int_0^T\|\partial u(t, \cdot)\|_{L^\infty}dt
\end{equation}
for classical solutions $u(t, x)$ and some $T > 0$. Cauchy-Swartz
inequality gives
\begin{equation}\label{g2}
\int_0^T\|\partial u(t, \cdot)\|_{L^\infty}dt \leq
T^{\frac{1}{2}}\Big(\int_0^T\|\partial u(t,
\cdot)\|_{L^\infty}^2dt\Big)^{\frac{1}{2}}.
\end{equation}
Note that the positive constant $T$ can be sufficiently small, we
only need to bound the quantity
\begin{equation}\label{g3}
\int_0^T\|\partial u(t, \cdot)\|_{L^\infty}^2dt.
\end{equation}
As in section 2, we prove the following Theorem:
\begin{thm}\label{thm71}
Suppose that the assumption $(\overline{H_1})$ in \eqref{g0} is
satisfied and $u(t, r)$ is the classical solution to the Cauchy
problem of the quasi-linear wave equation \eqref{a1}-\eqref{a2}.
Suppose furthermore that
\begin{eqnarray}\label{g4}
\frac{1}{c} \leq a(f) \leq c.
\end{eqnarray}
for some constant $c \geq 1$. Then there exists positive constants
$T
> 0$ and $C_j> 0$, $j = 7, 8, 9, 10$ depending only on $\big\|(f,
g)\big\|_{H^2(\mathbb{R}^3) \times H^1(\mathbb{R}^3)}$ such that
for $0 \leq t \leq T$, there holds

A priori estimate for $u$:
\begin{equation}\label{g5}
\|u(t, \cdot)\|_{L^\infty} \leq C_7.
\end{equation}

A priori estimates for $L_\pm v$:
\begin{equation}\label{g6}
\|L_+v(t, \cdot)\|_{L^\infty} \leq C_8,\quad  \|L_-v(t,
\cdot)\|_{L^\infty} \leq C_9.
\end{equation}

A priori estimate for $\partial u$
\begin{equation}\label{g7}
\int_0^T\|\partial u(t, \cdot)\|_{L^\infty}^2dt \leq C_{10}^2.
\end{equation}
\end{thm}

\begin{proof}
We use the same strategy as in the proof of Theorem \ref{thm21}.
By a straight forward calculation
\begin{eqnarray}\nonumber
|u(t, r) - f(r)| \leq \int_0^t|u_t(\tau, r)|d\tau \leq
C_{10}T^{\frac{1}{2}},
\end{eqnarray}
and noting that $T$ is a sufficiently small positive constant, we
deduce follows from \eqref{g4} and the smoothness of function
$a(\cdot)$ that
\begin{equation}\label{g8}
\begin{cases}
\frac{1}{2c} \leq a(u) \leq 2c,\\[-4mm]\\
\|u(t, \cdot)\|_{L^\infty} \leq 2\|f\|_{L^\infty} \leq
2\|f\|_{H^2(R^3)},
\end{cases}
\end{equation}
which gives the \textit{a priori} estimate for $u$ in \eqref{g5}.

By \eqref{d10}, we have
\begin{eqnarray}\nonumber
|L_+v(t, r)| \leq |L_+v(0, \beta)| + CC_9C_{10}T^{\frac{1}{2}},
\end{eqnarray}
which gives the \textit{a priori} bound depending only on
$\big\|(f, g)\big\|_{H^2(\mathbb{R}^3) \times H^1(\mathbb{R}^3)}$
for $L_+v$ if $T$ is sufficiently small. The similar manner gives
the \textit{a priori} bound for $L_-v$.

To get the \textit{a priori} estimate \eqref{g7}, we need the
following Lemma paralleling with Lemma \ref{lem36} and Lemma
\ref{lem37}:
\begin{lem}\label{lem71}
Suppose that \eqref{g4} and $(\overline{H_1})$ are satisfied. Then
\begin{equation}\label{g9}
\begin{cases}
\frac{r_-\big[\tau; \beta(t, r)\big] - r_-\big[\tau;
\beta\big(\alpha(t, r), 0\big)\big]}{r} \leq C,\\[-3mm]\\
\frac{1}{4} \leq \frac{dm(\tau, t)}{dt} \leq 4
\end{cases}
\end{equation}
holds for any $(t, r)$ with $t_+(r; 0) \leq t \leq T$ provided
that $T$ is sufficiently small, where $m(\tau, t)$ is given in
\eqref{c14}.
\end{lem}
We omit the proof of the above Lemma since it is the repeated
procedure as those of Lemma \ref{lem36} and Lemma \ref{lem37}.

By using Lemma \ref{lem71}, the \textit{a priori} estimate
\eqref{g7} can be proved by similar but simpler manner as in step
4 of section 6.
\end{proof}

Next, we show the desired stability of local strong solutions to
the Cauchy problem of the quasi-linear wave equation
\eqref{a1}-\eqref{a2}. Let $(f^i, g^i) \in H^2(\mathbb{R}^3)
\times H^1(\mathbb{R}^3)$, $i = 1, 2$ be two pairs of initial data
satisfying the assumption \eqref{g4} and $(\overline{H_1})$, and
$u_n^i(t, x)$, $i = 1, 2$ be the corresponding classical solutions
to the quasi-linear wave equation \eqref{a1} with
\begin{equation}\nonumber
u^i(0, x) = (\mathcal{J}_nf^i)(r), u_t^i(0, x) =
(\mathcal{J}_ng^i)(r),\quad x \in \mathbb{R}^3.
\end{equation}
To simplify our notations, we will drop the subscript $n$ as what
we have done in former sections. For $\lambda \in [0, 1]$, define
\begin{equation}\label{g10}
\begin{cases}
f^\lambda(x) = (1 - \lambda)f^1(r) + \lambda f^2(r),\\
g^\lambda(x) = (1 - \lambda)g^1(r) + \lambda g^2(r),
\end{cases}
\quad x \in \mathbb{R}^3.
\end{equation}
We further let $u^\lambda(t, r)$ be the unique global classical
solution to the quasi-linear wave equation \eqref{a1} with the
initial data
\begin{equation}\label{g11}
u^\lambda(0, x) = f^\lambda(r), u^\lambda_t(0, x) =
g^\lambda(r),\quad x \in \mathbb{R}^3.
\end{equation}

A straight forward calculations shows
\begin{eqnarray}\nonumber
&&\|\nabla u^1(t, \cdot) - \nabla u^2(t,
\cdot)\|_{L^2(\mathbb{R}^3)} +
  \|u_t^1(t, \cdot) - u_t^2(t,
  \cdot)\|_{L^2(\mathbb{R}^3)}\\\nonumber
&&= \Big\|\int_0^1\partial_\lambda\nabla u^\lambda(t,
  \cdot)d\lambda\Big\|_{L^2(\mathbb{R}^3)} +
  \Big\|\int_0^1\partial_\lambda u_t^\lambda(t,
  \cdot)d\lambda\Big\|_{L^2(\mathbb{R}^3)}\\\nonumber
&&\leq \sup_{0 \leq \lambda \leq 1}\big(\|
   \partial_\lambda\nabla u^\lambda(t,
  \cdot)\|_{L^2(\mathbb{R}^3)} + \|\partial_\lambda
  u_t^\lambda(t, \cdot)\|_{L^2
  (\mathbb{R}^3)}\big).
\end{eqnarray}
Thus, to show the stability result of the local strong solutions
to the quasi-linear wave equation \eqref{a1}, it suffices to prove
that
\begin{eqnarray}\label{g12}
\sup_{0 \leq \lambda \leq 1}\big(\|
   \partial_\lambda\nabla u^\lambda(t,
  \cdot)\|_{L^2(\mathbb{R}^3)} + \|\partial_\lambda
  u_t^\lambda(t, \cdot)\|_{L^2
  (\mathbb{R}^3)}\big) \leq  C\kappa
\end{eqnarray}
holds for $0 \leq t \leq T$, where
\begin{eqnarray}\label{g13}
\kappa = \big\|\nabla f^1 -
  \nabla f^2\big\|_{L^2(\mathbb{R}^3)} + \big\|g^1 -
  g^2\big\|_{L^2(\mathbb{R}^3)}.
\end{eqnarray}
Consequently, we turn to estimate the quantity
$\|\partial\partial_\lambda u^\lambda(t,
\cdot)\|_{L^2(\mathbb{R}^3)}$. This is the so-called homotopy
method.

First of all, similarly as in section 2, the standard energy
estimate gives that
\begin{eqnarray}\nonumber
&&\frac{d}{dt}E_1\big[\partial_\lambda u^\lambda(t)\big] \leq
C\|\partial u^\lambda(t, \cdot)\|_{L^\infty(\mathbb{R}^3)}
E_1\big[\partial_\lambda u^\lambda(t)\big]\\\nonumber &&\quad  +\
C\|\partial_\lambda u^\lambda(t, \cdot)\|_{L^\infty(\mathbb{R}^3)}
  E_1\big[\partial_\lambda u^\lambda(t)\big]^{\frac{1}{2}}
  E_2\big[u^\lambda(t)\big]^{\frac{1}{2}}.
\end{eqnarray}
Gronwall's inequality gives
\begin{eqnarray}\label{g14}
&&E_1\big[\partial_\lambda u^\lambda(t)\big]^{\frac{1}{2}}
  \\\nonumber
&& \leq \Big\{E_1\big[\partial_\lambda u^\lambda(0)
  \big]^{\frac{1}{2}} + \int_0^t \|\partial_\lambda u^\lambda(\tau,
  \cdot)\|_{L^\infty(\mathbb{R}^3)}E_2\big[u^\lambda(\tau)
  \big]^{\frac{1}{2}}d\tau\Big\}\\\nonumber
&&\quad \times \exp\Big\{\int_0^tC\|\partial
  u^\lambda(s, \cdot)\|_{L^\infty(\mathbb{R}^3)}ds\Big\}\\\nonumber
&&\leq C\Big\{\kappa + \int_0^t \|\partial_\lambda u^\lambda(\tau,
  \cdot)\|_{L^\infty(\mathbb{R}^3)}dt\Big\}\\\nonumber
&&\leq C\Big\{\kappa + t^{\frac{1}{2}}\Big(\int_0^T
  \|\partial_\lambda
  u^\lambda(\tau, \cdot)\|_{L^\infty(\mathbb{R}^3)}^2
  d\tau\Big)^{\frac{1}{2}}\Big\}.
\end{eqnarray}
It is clear that we need show
\begin{eqnarray}\label{g15}
\int_0^t\|\partial_\lambda u^\lambda(\tau,
  \cdot)\|_{L^\infty(\mathbb{R}^3)}^2d\tau \leq C\kappa^2
\end{eqnarray}
for $0 \leq t \leq T$. Note that $u^\lambda(t, r)$ satisfies the
estimates in Theorem \ref{thm71} since $f^\lambda$ and $g^\lambda$
satisfy the assumption \eqref{g4} and $(\overline{H_1})$. As the
proof for the estimate \eqref{g15} is similar as in step 4 of
section 6, we omit the details and just outline the ideas below.
We will drop the superscript $\lambda$ for notational convenience.

Differentiating the quasi-linear wave equation in \eqref{a1}, we
have
\begin{equation}\label{g16}
(r\partial_\lambda u)_{tt} - a(u)^2(r\partial_\lambda u)_{rr} =
2a(u)a^\prime(u)\partial_\lambda u(ru)_{rr}.
\end{equation}
Then rewrite the above equation as
\begin{equation}\label{g17}
M_-M_+v = \frac{a^\prime(u)}{a^2(u)} \Big[2a(u)\partial_\lambda
u(ru)_{rr} - rM_-u\partial_\lambda u_t\Big].
\end{equation}
For $t \geq t_+(r; 0)$, integrating the above equation \eqref{g17}
along the corresponding plus and minus characteristics gives
\begin{eqnarray}\label{g18}
&&\partial_\lambda u(t, r) =
\frac{1}{r}\int_0^rM_+(r\partial_\lambda u)\big[0,
\beta\big(t_+(\sigma;
  \alpha), \sigma\big)\big]d\sigma - \frac{1}{r}\int_0^r
  \int_0^{t_+(\sigma; \alpha)}\\\nonumber
&&\quad \frac{a^\prime(u)}{a(u)} \Big[2a(u)\partial_\lambda
  u(ru)_{rr} - rM_-u\partial_\lambda u_t\Big]\Big(s, r_-\big[s;
  \beta\big(t_+(\sigma; \alpha), \sigma\big)\big]\Big)dsd\sigma,
\end{eqnarray}
For $t < t_+(r; 0)$, integrating the equation \eqref{g17} along
the corresponding plus and minus characteristics gives
\begin{eqnarray}\label{g19}
&&\partial_\lambda u(t, r) =
\frac{1}{r}\int_0^rM_+(r\partial_\lambda u)\big[0,
\beta\big(t_+(\sigma;
  \gamma), \sigma\big)\big]d\sigma - \frac{1}{r}\int_0^r
  \int_0^{t_+(\sigma; \gamma)}\\\nonumber
&&\quad \frac{a^\prime(u)}{a(u)} \Big[2\partial_\lambda
  u(ru)_{rr} - rM_-u\partial_\lambda u_t\Big]\Big(s, r_-\big[s;
  \beta\big(t_+(\sigma; \gamma), \sigma\big)\big]\Big)dsd\sigma,
\end{eqnarray}
Following the similar manner as in step 4 of section 6, we can
complete the proof of the estimate \eqref{g15}.

\bibliographystyle{amsplain}

\end{document}